\documentclass[11pt,leqno]{article}
\usepackage{amssymb}
\usepackage[monochrome]{color}
\newtheorem{proposition}{Proposition}[section]
\newtheorem{lemma}{Lemma}[section]
\newtheorem{corollary}{Corollary}[section]
\newtheorem{remark}{Remark}[section]

\newcommand{\dfrac}{\displaystyle\frac}
\newcommand{\dint}{\displaystyle \int}
\renewcommand{\H}{{\cal H}}
\newcommand{\T}{{\cal T}}
\newcommand{\D}{{\cal D}}
\newcommand{\W}{{\cal W}}
\newcommand{\K}{{\cal K}}
\renewcommand{\L}{{\Lambda}}
\newcommand{\ve}{{\varepsilon}}

\newcommand{\J}{{\cal J}}
\definecolor{darkred}{rgb}{0.8,0,0}

\newcommand{\qed}{\hfill$\square$\vspace{0.3cm}}
\begin{document}
\title{\bf Functional aspects of the Hardy inequality. Appearance of a hidden energy}
\author{\Large J. L. V\'azquez, Madrid\footnote{juanluis.vazquez@uam.es} \\[6pt]
\Large N. B. Zographopoulos, Athens
\footnote{zographopoulosn@sse.gr, nzograp@gmail.com}
}
\date{}
\maketitle
\pagestyle{myheadings} \thispagestyle{plain} \markboth{}{}
\maketitle
\begin{abstract}

Starting with a functional difficulty appeared in the paper \cite{vz00} by V\'azquez and Zuazua, we obtain new insights  into the Hardy Inequality and the evolution problem associated to it by
means of a reformulation of the problem. Surprisingly, the
connection of the energy of the new formulation with the standard
Hardy functional is nontrivial, due to the presence of a Hardy singularity energy.
This corresponds to a loss for the total energy. The problem arises when the equation is posed in a bounded domain, and also when posed in the whole space.

We also  consider an equivalent problem with inverse square potential on an exterior domain. The extra energy term is then present as an effect that comes from infinity, a kind of hidden energy. In this case, in an unexpected way, this term is  additive to the total energy, and it may even constitute the main part of it.
\end{abstract}
\section{Introduction}
\setcounter{equation}{0}
In this paper we contribute some new results on the Hardy
Inequality posed in a  bounded domain, in an exterior domain of
$\mathbb{R}^N$, or the whole $\mathbb{R}^N$ and on the corresponding parabolic evolution. The
motivation came from a functional difficulty we found in the work
\cite{vz00}, where the following singular evolution problem  was
studied:
\begin{equation}\label{eq1}
\left\{\begin{array}{ccc}
u_t &=& \Delta u + c_*\, \dfrac{u}{|x|^2},\; x \in \Omega,\; t>0,\\
u(x,0) &=& u_0(x),\;\; \mbox{for}\;\; x \in \Omega, \nonumber \\ u(x,t)&=&0\;\;\mbox{in}\;\;\partial\Omega,\;t>0\,.
\end{array}\right.
\end{equation}
with critical coefficient $c_*=(N-2)^2/4$. The space dimension is $N\ge 3$ and $\Omega$ is a bounded domain in $\mathbb{R}^N$ containing $0$, or $\Omega = \mathbb{R}^N$.
More precisely, the authors in \cite{vz00} studied the well-posedness and described the asymptotic
behavior of (\ref{eq1}). Moreover, they obtained improved Hardy inequalities and completed the study of
the spectrum of the associated  eigenvalue problem. This problem is
closely connected with the Hardy inequality:
\begin{equation} \label{hardyineq}
\int_{\Omega} |\nabla u|^2\, dx > \biggr ( \frac{N-2}{2} \biggr
)^2 \int_{\Omega} \frac{u^2}{|x|^2}\, dx,
\end{equation}
which is well known to hold for any $\phi \in
C_0^{\infty}(\Omega)$. For Hardy type inequalities and related
topics we refer to \cite{bm, d99, kmp, maz85, ok}. Due to this connection,
$c_*$, which is the best constant in the inequality, is  also
critical for the basic theory of the evolution equation. Indeed,
the usual variational theory applies to the subcritical cases: \
$u_t = \Delta u +  c\, u/|x|^2$ \ with $c<c_*$, using the standard
space $H^1_0(\Omega)$, and a global in time solution is then
produced. On the other hand, there are no positive solutions of
the equation for $c>c_*$ (instantaneous blow-up),
\cite{BG,cm99}. In the critical case we still get existence
but the functional framework changes; this case serves as an
example of interesting functional analysis and more complex
evolution.

In order to analyze the behavior of the solutions of Problem
(\ref{eq1}) in \cite{vz00}, the Hardy functional
\begin{equation} \label{Hardyfunctional}
I_\Omega[\phi] := \int_{\Omega} |\nabla \phi|^2\, dx -
\left ( \frac{N-2}{2} \right )^2 \int_{\Omega}
\frac{\phi^2}{|x|^2}\, dx\,,
\end{equation}
is considered as the Dirichlet form naturally  associated to the
equation. This form is positive and different lower bounds have
been obtained see \cite{bv97, vz00}. Note
that the expression is finite for $u\in H^1_0(\Omega)$, but it can
also be finite as an improper integral for other functions having
a strong singularity at $x=0$, due to cancelations between the
two terms. To take this possibility into account, the  Hilbert
space $H$ was introduced in \cite{vz00} as  the completion of the
$C_{0}^{\infty} (\Omega)$ functions under the norm
\begin{equation} \label{oldH}
||\phi||^{2}_{H(\Omega)} = I_{\Omega} [\phi],\;\;\; \phi \in C_{0}^{\infty} (\Omega).
\end{equation}
According to Section 5 of \cite{vz00}, this space allows us to
define in a natural way  a self-adjoint extension of the
differential operator $L(u):= -\Delta u- c_* \, u/|x|^2$
(the Friedrichs extension) and then to
use standard theory to generate a semigroup and describe the
solutions using the spectral analysis. The study of the spectrum
leads to an associated elliptic eigenvalue problem, the solution
of which turns out to be a classical problem in separation of
variables.

\medskip

\noindent {\sc Problem with the singularities.} The  separation of
variables analysis produces some singular solutions. In
particular, the maximal singularity (corresponding to the first
mode of separation of variables) behaves like $|x|^{-(N-2)/2}$
near $x=0$, and this function is not in $H^1_0(\Omega)$. Now, this
solution must belong to the space $H$ associated to the quadratic
form, hence the conclusion $H\ne H^1_0(\Omega)$.   We recall that
this is a peculiar phenomenon of the equation with critical
exponent $c_*=(N-2)^2/4$. For values of $c< c_*$ the maximal
singularity is still in $H^1_0(\Omega)$.

However, there
must be a gap in the argument of \cite{vz00}. We have realized
that with the proposed definition of $H$, there exists a problem
with the solutions of the evolution problem  having the maximal
singularity. The verification is quite simple in the case where
$\Omega =B_1$, the unit ball in $\mathbb{R}^N$ centered at the
origin. Then, the minimization problem
\begin{equation}\label{01}
\min_{u \in H} \frac{||u||^2_{H}}{||u||^2_{L^2}}
\end{equation}
has as a solution the function
\begin{equation}\label{e1}
e_1 (r) = r^{-(N-2)/2}\, J_0 (z_{0,1}\, r),\;\;\; r=|x|,
\end{equation}
$J_0$ is the Bessel function with $J_0 (0) =1$, up to normalization and $z_{0,1}$ denotes the first zero of $J_0$. This function plays a big role in the asymptotic behavior of general solutions of Problem (\ref{eq1}). The minimum value of (\ref{01}) is
\[
\mu_1 = z^2_{0,1}.
\]
Moreover, the quantity $I_{B_1}(e_1)$ is well defined as a principal value. Assuming that
\begin{equation}\label{H=I}
||e_1||^2_H = I_{B_1}(e_1),
\end{equation}
from the definition of $H$, for any $\varepsilon>0$, we should find a $C_0^{\infty}$-function $\phi$,
such that
\begin{equation}\label{02}
||e_1-\phi||^2_{H} < \varepsilon.
\end{equation}
However, setting
\[
e_1(r) - \phi(r) = |r|^{-(N-2)/2} w(r),\qquad r=|x|\,,
\]
we have
\[
w(r) = J_0 (z_{0,1}\, r) - r^{(N-2)/2} \phi(r)\,,
\]
which is regular at $r=0$, and after some straightforward
calculations we find that, using as norm the square
root of $I_{B_1}$,
\begin{equation} \label{04}
||e_1-\phi||^2_{H} = \int_{\Omega} |x|^{-(N-2)} |\nabla w|^2\, dx +
\frac{1}{2} \int_{\Omega} \nabla |x|^{-(N-2)} \cdot \nabla w^2\, dx.
\end{equation}
The last integral is not zero, due to the presence of a "boundary term",
when we integrate by parts. More precisely,
\begin{equation} \label{04a}
\frac{1}{2} \int_{\Omega} \nabla |x|^{-(N-2)} \cdot \nabla w^2\, dx =
\frac{N(N-2)}{2}\, \omega_N\, J^2_0(0) =
\frac{N(N-2)}{2}\, \omega_N.
\end{equation}
where $\omega_N$ denotes the Lebesgue measure of the unit ball in
$\mathbb{R}^N$. Thus,
\begin{equation}\label{03}
||e_1-\phi||^2_{H} \geq \frac{N(N-2)}{2}\, \omega_N,
\end{equation}
which contradicts (\ref{02}). Thus, we see that $e_1$ fails to be in $H$, since it cannot be approximated by $C_0^{\infty}$-functions and this will  happen for every function with the maximal singularity.

Therefore, under the assumption (\ref{H=I}), the space $H$ seems not to be correctly defined in \cite{vz00} to apply the rest of the theory, since there exists a problem in dealing with very singular behavior near $x=0$ that is not covered by approximation with infinitely smooth functions. Actually, this was our first impression.

\medskip

\subsection*{New results}

\noindent {\bf 1}.  The examination of the difficulty shows that
the proposed norm $I_\Omega$ is too detailed near the singularity
and produces a topology that is too fine to allow the convergence
of $\phi_n\in C^\infty_c(\Omega)$ to $e_1$. By means of a
transformation already proposed in \cite{bv97}, we obtain a more suitable
norm $N$ that is equivalent to $I_\Omega^{1/2}$ on $C_0(\Omega)$,
but is gross enough near the singular point. In this way, we are
able to define a possibly larger closure, that we call $\cal H$,
that contains all the functions needed for constructing the
evolution. This is done in Subsections \ref{sub.definition} and
\ref{sec.sfp}.

We proceed next to re-examine the above mentioned difficulty. We
will show that the spaces $H$ and $\cal H$ are indeed the same.
What is different is the norm that was implicitly assumed to be
acting in $H$ for solutions that do not necessarily vanish at
$x=0$, which in principle seemed to be $I_{\Omega}^{1/2}$ taken in
the sense of principal value. When both terms of $I_\Omega$ become
infinite the correct definition of the norm is a particular limit
that we call the {\sl cutoff limit}. This is explained in
Subsections \ref{sec.hf} - \ref{sub.normH} where we examine the
connection of the new norm with the Hardy functional; the
difference is characterized in terms of a certain value the {\sl
Hardy singularity energy} (HS energy for short) that we precisely define. We think
that the existence of the two different norms that coincide on
$C_c(\Omega\setminus \{ 0\})$ is quite interesting and was
unexpected for us.

The spectral analysis of \cite{vz00} becomes rigorous in this setting, and a contraction semigroup is associated to the evolution problem by standard variational methods in weighted Sobolev spaces, up to unitary equivalence. This analysis is carefully explained in Section \ref{sec.evol}.

\medskip

{\bf 2}. In Section \ref{sec.h} we discuss a result which has its
own interest; the Critical Caffarelli--Kohn--Nirenberg
Inequalities, in a bounded domain. It was well known (see
\cite{cw01}) that these inequalities are related to the Hardy
inequality with $c < c_*$. The critical case is as expected
related with the Hardy inequality with $c = c_*$. The
proper functional setting that we had to consider for Problem
(\ref{eq1}) leads us naturally to these critical inequalities. We
also give the connection of this new space with the Sobolev space
$D^{1,2} (\mathbb{R}^N)$, with the use of a proper transformation.
Based on this transformation, we describe an easy to apply
argument concerning the existence of non-$H_0^1$ minimizers,
establishing that their behavior at the origin, is precisely
$|x|^{-(N-2)/2}$.

\medskip

{\bf 3}. In Section \ref{sec.exd} we explore the existence of an analogue  of the
Hardy singularity energy for problems posed in
exterior domains. The Kelvin transform suggests that the most natural problem to study is the following:
\begin{equation}\label{exd1}
\left\{\begin{array}{ccll}
|y|^{-4}\, w_t (y,t) &=& \Delta w (y,t) + c_*\, \dfrac{w (y,t)}{|y|^2},\ &y \in B^c_{\delta},\; t>0,\\
w(y,0) &=& w_0(y),\;\; &\mbox{for}\;\; y \in B^c_{\delta}, \nonumber \\
w(y,t)&=& 0\;\; &\mbox{for} \;\;
|y|=\delta,\;t>0\,.
\end{array}\right.
\end{equation}
where $c_*=(N-2)^2/4$ is the critical coefficient,  $B^c_{\delta} = \mathbb{R}^N \backslash B_{\delta} (0)$ is the standard exterior domain and $\delta>0$. Without loss of generality we take $\delta=1$.

Problem (\ref{exd1}) has the  striking property  that the Hardy functional posed in the exterior domain is not necessarily  a positive quantity; we will show that for functions which vanish at $\partial B^c_{\delta} (0)$ and behaving at infinity like $|y|^{-(N-2)/2}$ it may be negative. We avoid the difficulty by basing
our existence theory  on the unitary equivalence via the Kelvin transform. Hardy inequalities on exterior domains were studied in the work \cite{bft}.

The novel feature of the bounded domain, namely, the HS energy term, does exist also in the case of the exterior domain, but it appears  at infinity. Besides, there is a big difference with the bounded domain case since in this case the new energy is not only additive to the total energy involved in the evolution; it may even represent the main part of it. The energy term at infinity looks to us like a ``hidden'' or ``dark" energy. We are not in a position to make a physical interpretation, but we have proved the existence of a such energy and give its precise formula.  This seems to be the first study of an evolution problem with such curious properties in an exterior domain.

Equations (\ref{eq1}) and (\ref{exd1}) are linear parabolic equations. In this sense, it is expected that their dynamics should be trivial. However, the presence of a singular potential with an inverse square power singularity changes things; due to that, the first equation has interesting behavior  at zero and the second at infinity. Moreover, we treat the limiting case of the best constant $c_*$, in the sense of the Hardy inequality. The singularity and the best constant make the dynamic of these equations far from being trivial. The appearance of correcting terms in the form of hidden energies
has been unexpected to us. We think that the question deserves further investigation.\normalcolor

\medskip

\noindent {\bf 4.} In Section \ref{sec.rn} we consider the Hardy functional and the corresponding evolution problem posed on $\mathbb{R}^N$.  The problem was  studied in \cite{vz00}, see also \cite{vzog}.
It is known that the Hardy inequality on $\mathbb{R}^N$ is sharp, i.\,e., it cannot be improved, at least up to some $L^p$-norm. This may be also seen as consequence of the transformation $u (x)= |x|^{-(N-2)/2}\, v(x)$; we can find functions, not belonging to $L^2 (\mathbb{R}^N)$, such that the Hardy inequality holds. To overcome this difficulty, the authors in
\cite{vz00} made use of the similarity variables. In this case we do not find any HS energy at infinity, as in the previous exterior domain model, see \cite{vzog}.

The main result of the section is to prove that the Hardy functional on $\mathbb{R}^N$, for a certain class of functions, maybe improved by the $L^2 (\mathbb{R}^N)$-norm, i.e. a Hardy-Poincar\'{e} inequality on $\mathbb{R}^N$, holds in this sense.  The idea is to use a new, proper, transformation.
\begin{lemma}  \label{lemHPRN} \textbf{\emph{(Hardy-Poincar\'{e} inequality on $\mathbb{R}^N$)}}
Let $v$ be any function in $C^{\infty}_0 (\mathbb{R}^N)$, $N \geq
3$ and let
\begin{equation} \label{trans1}
u = |x|^{-\frac{N-2}{2}}\, J_0 (|x|)\, v,
\end{equation}
where $J_0$ is the Bessel function with $J_0 (0) =1$, up to normalization. Then,
\begin{equation} \label{HardyPoincare}
\int_{\mathbb{R}^N} |\nabla u|^2\, dx - \biggr ( \frac{N-2}{2} \biggr
)^2 \int_{\mathbb{R}^N} \frac{u^2}{|x|^2}\, dx > \int_{\mathbb{R}^N} u^2\, dx.
\end{equation}
The best constant in {\rm (\ref{HardyPoincare})} is $1$ and there exists no minimizer.
\end{lemma}

As a consequence of the above Lemma we derive a result, see Remark
\ref{remR2}, which has its own interest. The connection of
transformation $u (x)= |x|^{-(N-2)/2}\, v(x)$ with the reduction
of dimension from $N$ to $2$ is one of the leitmotifs of paper
\cite{bv97} and will be used frequently
in the sequel. This result concerns the improvement, for a certain
class of functions, of the ``norm" of $D^{1,2} (\mathbb{R}^2)$.
Note that $D^{1,2} (\mathbb{R}^2)$ is not a well defined space.

Lemma \ref{lemHPRN} is, as far as we know, the first result
concerning improvement, with a norm, of the Hardy inequality on
$\mathbb{R}^N$. The only result, up to our knowledge, that gave an
improvement of this inequality is in the work \cite{cf08}, where
the authors obtained a non-standard improvement.

Returning now to the Hardy functional, we define a new weighted
space $\J$ through (\ref{trans1}). We state some properties of
this space, focused on the behavior at infinity; the functions
belonging in this space belong also to $L^2 (\mathbb{R}^N)$, so
they vanish at infinity. It turns out that $\J$ is a well defined
Hilbert space. The connection of this space with the Hardy
functional on $\mathbb{R}^N$ is nontrivial, as it contains the
Hardy's singularity energy at the origin, but it also contains
singularity energy at each zero $z_m$ of the Bessel function.
Transformation (\ref{trans1}) removes the singularity at infinity.
However, the cost we have to pay, is the existence of an infinite
series of singularities at $z_m$. Contrary to the case of the
origin, the singularity energy at these points, increases the
energy.

Some further results on Hardy type inequalities and the related hidden energies are contained in the forthcoming work \cite{vzog}.

\section{Proper functional setting. Bounded domain case}\label{sec.pfs}
\setcounter{equation}{0}

We start the detailed study  by analyzing the case of a bounded
domain $\Omega\subset \mathbb{R}^N$, $N \geq 3$.


\subsection{Transformation and definition of spaces} \label{sub.definition}

The way we follow to address the difficulty mentioned in the
introduction and to properly pose Problem (\ref{eq1}) is to
introduce a more convenient variable by means of the formula
\begin{equation} \label{trans}
u (x)= |x|^{-(N-2)/2}\, v(x).
\end{equation}
We will write the transformation as $u=\T(v)$. Clearly, this is an
isometry from the space \ $X=L^2(\Omega)$ \ into the space \
$\widetilde X=L^2(d\mu,\Omega)$, $d\mu=|x|^{2-N}dx$. Many
arguments of \cite{vz00} were also based on  transformation
(\ref{trans}), which was first used in \cite{bv97} and then in
many papers concerning results on Hardy's inequalities. The great
advantage of this formula is that it simplifies $I_{\Omega}(u) $,
 at least for smooth functions, into
\begin{equation}\label{qf.v}
I_1(v) :=\int_{\Omega} |x|^{-(N-2)} |\nabla v|^2\, dx\,.
\end{equation}
It is easy checked that $I_{\Omega}(u) = I_1(v)$ for functions $u
\in C_{0}^{\infty} (\Omega)$. However, the equivalence fails for
functions with a singularity like $|x|^{-(N-2)/2}$ at the origin,
as we have hinted before and will explain below in detail. Our
proposal is to use this formulation for the definition of the new
space, $\H$. An important observation is that when $u(x,t)$ is a
solution of equation (\ref{eq1}), then $v$ satisfies the following
associated equation
\begin{equation} \label{Cauchy-Dirichlet_v}
v_t = |x|^{N-2}\; \nabla \cdot \left( |x|^{-(N-2)}\, \nabla v
\right),
\end{equation}
with clear equivalence for $x \ne 0$. This last form gives the
clue to the proper variational formulation to be followed here.
First, the space associated to this equation through the quadratic
form (\ref{qf.v}), is defined as the weighted space \ $\widetilde
{\cal H} = W^{1,2}_0(d\mu,\Omega)$, which is the
completion of the $C_0^{\infty}(\Omega)$ functions under the norm
\begin{equation} \label{weigted space}
||v||^2_{\widetilde {\cal H}} = \int_{\Omega} |x|^{-(N-2)}\,
|\nabla v|^2\, dx\,.
\end{equation}
Following the usual procedure of the Calculus of Variations, we
take an appropriate  base space  which is
$\widetilde X=L^2(d\mu,\Omega)$, and then the quadratic form
(\ref{qf.v}) has as form domain the subspace $\widetilde {\cal H}$
where $I_1(v)$ is finite. Then, it can be proved  that $L(v):=-
|x|^{N-2}\; \nabla \cdot \left( |x|^{-(N-2)} \nabla v\right)$ is a
positive self-adjoint operator in the space $D(L)=\{v\in
\widetilde {\cal H}: L(v)\in \widetilde X \}$. It is also known that
$D(L^{1/2})=\widetilde {\cal H}$. Hence, the variational approach
works for $v$. See further analysis below.

We translate these results to the original framework. $\H$ is
defined as the isometric space of $\widetilde {\cal
H}=W^{1,2}_0(|x|^{-(N-2)}dx,\Omega)$ under the transformation $\T$
given by (\ref{trans}). In other words, $\H$ is defined as the
completion of the set
\[
\left \{ u = |x|^{-\frac{N-2}{2}}\, v,\;\;\; v \in C_{0}^{\infty}
(\Omega) \right \}=\T (C_{0}^{\infty}(\Omega)),
\]
under the norm $N(u)=\|u\|_{\cal H}$ defined by
\begin{equation} \label{newH}
||u||^2_{\H} = \int_{\Omega}
|x|^{-(N-2)}\, |\nabla (|x|^{\frac{N-2}{2}}\, u)|^2 \, dx.
\end{equation}
%
%
%
\subsection{On the new space}\label{sec.sfp}

From the definition it follows that  $\H \supset H_0^1 (\Omega)$.
To show that $H_0^1 (\Omega)\ne \H$, we consider any  function in $\cal
H$ such that $u(x) \sim |x|^{-\frac{N-2}{2}}$ at the origin; it cannot belong to $H^1_0(\Omega)$. We continue with two
rather simple observations before examining the important question
of relation of the present approach with the old one, where we
will obtain a striking correction term.

\medskip

\noindent $\bullet$ {\sl  The functions that behave at the origin
like $|x|^{-(N-2)/2}$ are not the most singular in $\H$. } This
was mentioned in passing in \cite{vz00}. We give a counterexample
here. Consider a function $w$ which behave near the origin like
$$
v(x)\sim (\log 1/|x|)^{a}, \qquad 0<a<1/2.
$$
It is easy to check that $v \in \widetilde {\cal H}$. Therefore,
$|x|^{-(N-2)/2}\, v$ belongs to $\H(\Omega)$, even if
$v(x)\to\infty$ as $x\to 0$.

\medskip

\noindent $\bullet$ {\sl There exist functions that belong in
$\bigcap_{q<2} W^{1,q}(\Omega)$ but they do not belong to $\H$.}
For example, a function $v$ that behaves near the origin like
\[
v(x) \sim (\log (1/r))^{1/2} ,
\]
fails to be in $\widetilde {\cal H}$, hence
\[
u \sim r^{-\frac{N-2}{2}} (\log (1/r))^{1/2} ,
\]
does not belong to $\H$. On the other hand, $u \in
W^{1,q}(\Omega)$, for all \  $1 \leq q
<2$.

%
\subsection{Connection of space $\H$ with the Hardy functional} \label{sec.hf}

By Hardy functional we refer to $I_{\Omega}(u)$ defined in
(\ref{Hardyfunctional}) with the integral defined in
the sense of principal value around the origin when both separate
integrals diverge. Denote by $B_\varepsilon$ the ball
centered at the origin with radius $\varepsilon$ and by
$B^c_\varepsilon$ its complement in $\Omega$. Assume now that $u
\in \H$, so that $v= |x|^{(N-2)/2}\, u \in \widetilde {\cal H}$.
Then, we have that
\begin{equation} \label{I_varepsilon}
I_{B^c_\varepsilon} [u] = \int_{B^c_\varepsilon} |\nabla u|^2\, dx
- \left ( \frac{N-2}{2} \right )^2 \int_{B^c_\varepsilon}
\frac{u^2}{|x|^2}\, dx.
\end{equation}
By change of variables and integration by parts the following
remarkable formula is obtained:
\begin{equation}\label{twoint}
I_{B^c_\varepsilon} [u] = \int_{B^c_\varepsilon}
|\nabla v|^2\, |x|^{2-N} dx+ \frac{N-2}{2}\,
\varepsilon^{-(N-1)}\, \int_{S_\varepsilon} v^2\, dS\,,
\end{equation}
where $dS$ is the surface measure. Next, we denote by
$\L_{\varepsilon}$ the quantity:
\begin{eqnarray} \label{Lepsilon}
\L_{\varepsilon} (u) &=& \frac{N-2}{2}\, \varepsilon^{-(N-1)}\,
\int_{S_\varepsilon} v^2\, dS\, \nonumber \\ &=& \frac{N-2}{2}\, \varepsilon^{-1}\,
\int_{S_\varepsilon} u^2\, dS\,,
\end{eqnarray}
that represents a kind of {\sl Hardy energy at the singularity}. Is clear
that
\[
\lim_{\varepsilon \to 0} \int_{B^c_\varepsilon} |\nabla v|^2\,
|x|^{2-N} dx = ||v||^2_{\widetilde {\cal H}}.
\]
In order to take the limit $\varepsilon \to 0$, in
(\ref{I_varepsilon}) we distinguish the following cases:

\medskip

\noindent $\bullet $ If $u \in H_0^1 (\Omega)$, then $u\in {\cal
H}$ and we have
\[
\L(u):= \lim_{\varepsilon \to 0} \L_{\varepsilon} (u) =0,
\]
thus the limit as $\varepsilon \to 0$, in (\ref{I_varepsilon}),
gives the well known formula
\[
I_\Omega [u] = ||v||^2_{\widetilde {\cal H}}:=N^2(u),
\]
which holds for any $u \in H_0^1 (\Omega)$.  {\sl Note that the
converse is not true;} If $\L(u)=0$, it does not imply that $u \in
H_0^1 (\Omega)$. For example, take a function $u$ such that $v$
behaves at zero like $(-\log |x|)^{-1/2}$.

\medskip

\noindent $\bullet $ If $v \in \widetilde {\cal H}$ is such that
$\lim_{|x| \to 0} v^2(x) = v^2 (0)$ exists as a real positive
number; then $u\in \H$ but $u \not\in H_0^1 (\Omega)$. In this
case
\begin{equation} \label{Lfinite}
\L(u) = \frac{N(N-2)}{2}\, \omega_N\, v^2(0),
\end{equation}
where $\omega_N$ denotes the Lebesgue measure of the unit ball in
$\mathbb{R}^N$.  $\L(u)$ is then a well defined positive number
and (\ref{I_varepsilon}) implies that
\begin{equation} \label{Ifinite}
I_\Omega [u] = ||v||^2_{\widetilde {\cal H}} + \L(u).
\end{equation}
We note that this is the case of the principal eigenfunction of
(\ref{eigen_u}) and the case of the minimizer of the Improved
Hardy-Sobolev inequality, see \cite{z10}, in the radial case.
As it will be clear in Subsection \ref{sub.furtherprop},
this is the case for the minimizers of
\begin{equation} \label{otherminimizers}
   \min_{u \in H} \frac{||u||_{H}}{||u||_{L^p}},\;\;\;\;\;\; 1 \leq p < \frac{2N}{N-2}.
\end{equation}
Note also that in this case the inner product of $\H$ is given by
\begin{eqnarray}\label{innerproduct}
<u_1,u_2>_{\H} &=& \int_{\Omega} \nabla u_1 \cdot \nabla u_2\, dx -
\left ( \frac{N-2}{2} \right )^2 \int_{\Omega} \frac{u_1\,
u_2}{|x|^2}\, dx \nonumber \\
&&- \frac{N(N-2)}{2}\, \omega_N\, v_1(0)\, v_2(0).
\end{eqnarray}

\medskip

\noindent $\bullet $ If $v \in \widetilde {\cal H}$ is such that
$v$ at zero is bounded but the $\lim_{x \to 0} v^2(x)$ does not
exist, i.\,e., $v$ oscillates near zero. For example, let
\[
v \sim \sin \left( (-\log |x|)^a \right ),\;\;\; |x| \to 0.
\]
Then, $v$ belongs in $\widetilde {\cal H}$ for some $0<a<1/2$, so
$u = |x|^{-(N-2)} v \in \H$. In this case, the limit $L(u)$ does
not exist, since it oscillates, and from (\ref{I_varepsilon}) we
have that the same happens to the Hardy functional, in the sense
that
\begin{equation} \label{limitHardy}
\lim_{\varepsilon \to 0} \left ( I_{B^c_\varepsilon} [u] -
\L_{\varepsilon} (u) \right ) = ||v||^2_{\widetilde {\cal H}}.
\end{equation}

\noindent $\bullet $  If $v \in \widetilde {\cal H}$ is such that
$\lim_{x \to 0} v^2(x) = \infty$. For example, let
\[
v \sim (-\log |x|)^a ,\;\;\; |x| \to 0.
\]
Then, $v$ belongs in $\widetilde {\cal H}$ for some $0<a<1/2$, so
$u = |x|^{-(N-2)} v \in \H$. Is clear that, $\L(u) = \infty$, and from
(\ref{I_varepsilon}) we have that the same happens to the Hardy
functional, in the sense that (\ref{limitHardy}) holds.

\medskip

Note that in all cases $\L_{\varepsilon}$ is a nonnegative
quantity, for every $\varepsilon >0$ and so is
$I_{B^c_\varepsilon} [u]$. As a consequence, we obtain a
generalized form of the Hardy inequality valid in the limiting
case of (\ref{limitHardy}), when the Hardy functional is not
defined or is infinite. We do not know if there is any physical
meaning for the singularity energy we have found. It looks like an
energy defect at the singularity.

\subsection{The spaces $\H$ and $H$ are the same} \label{sub.H=H}
We recall that $H$ was introduced as  the completion of the
$C_{0}^{\infty}(\Omega)$ functions under the norm
$I_\Omega^{1/2}$. The proof of $\H=H$ relies on showing that the
set $C_0^{\infty} (\Omega \backslash \{0\})$ is a dense set in
both spaces, and on observing that the two norms coincide on that
subset.

The following lemma follows from \cite[Lemma 2.1]{cw01}, which
studies the subcritical case.
\begin{lemma} \label{denseH}
Holds that $C_0^{\infty} (\Omega \backslash \{0\})$ is a dense set in $H$.
\end{lemma}
\emph{Proof} By the definition of $H$, it suffices to prove that
\[
C_0^{\infty} (\Omega) \subset \overline{C_0^{\infty} (\Omega \backslash \{0\})}^{||.||_H}
\]
Let $\rho(t)$ be a cutoff function that is $1$ for $t \geq 2$ and
$0$ for $0<t\leq 1$, $0\le \rho(t)\le 1$. For a fixed $u \in
C_0^{\infty} (\Omega)$ we define
\begin{equation} \label{uepsilon}
u_{\varepsilon}(x) = \rho (|x|/\varepsilon)\, u(x) \in C_0^{\infty} (\Omega \backslash \{0\}).
\end{equation}
Then we have that
\begin{eqnarray*}
\left ( \int_{\Omega} \left ( \nabla u_{\varepsilon} - \nabla u \right )^2\, dx \right )^{1/2} &=& \left (  \int_{\Omega} \left ( \nabla \left (\rho (|x|/\varepsilon) \right )\, u + (\rho-1) \nabla u \right )^2\, dx \right )^{1/2} \\ &\leq& \left (  \int_{A_{\varepsilon}} \left ( \nabla \rho (|x|/\varepsilon) \right )^2\, u^2\, dx \right )^{1/2} + \left (  \int_{B_{2\varepsilon}} (\rho-1)^2 |\nabla u|^2\, dx \right )^{1/2},
\end{eqnarray*}
where $A_{\varepsilon} = \{x \in \Omega, \varepsilon < |x| < 2\varepsilon  \}$. The first integral is estimated as
\[
\int_{A_{\varepsilon}} \left ( \nabla \rho (|x|/\varepsilon) \right )^2\, u^2\, dx \leq c\, ||u||^2_{L^{\infty}}\, \varepsilon^{-2}\, \int_{A_{\varepsilon}} (\rho')^2\, dx \leq c\, \varepsilon^{N-3}\, \int_{\varepsilon}^{2\varepsilon} (\rho')^2\, dr \to 0,
\]
as $\varepsilon \downarrow 0$. Also, the second integral tends to zero, since
\[
\int_{B_{2\varepsilon}} (\rho-1)^2 |\nabla u|^2\, dx \leq c\, \int_{B_{2\varepsilon}} dx,
\]
so
\begin{equation} \label{H1}
 \int_{\Omega} \left ( \nabla u_{\varepsilon} - \nabla u \right )^2\, dx \to 0,
\end{equation}
as $\varepsilon \downarrow 0$. In addition,
\[
\int_{\Omega} \frac{(u_{\varepsilon}-u)^2}{|x|^2}\, dx \leq c\, \varepsilon^{N-3}\, \int_{B_{\varepsilon}} dx \to 0,
\]
as $\varepsilon \downarrow 0$. Finally, we conclude that
\[
||u_{\varepsilon}-u||_H \to 0,
\]
and the proof is completed. $\blacksquare$ \vspace{0.3cm}

Next we will prove that the $C_0^{\infty} (\Omega \backslash \{0\})$-functions are also dense in $\widetilde
{\cal H}$ and hence in $\H$. However, the functions defined by (\ref{uepsilon}) are not good approximations in $\widetilde {\cal H}$.
In other words, the arguments of \cite[Lemma 2.1]{cw01} cannot hold in the critical case $a=(N-2)/2$.
\begin{lemma} \label{nondense}
Let $v \in C_0^{\infty} (\Omega)$, $v(0) \ne 0$. Then the $v_{\varepsilon}$, given by $v_{\varepsilon}(x) = \rho (|x|/\varepsilon)\, v(x)$, cannot approach $v$ in  $\widetilde {\cal H}$.
\end{lemma}
\emph{Proof} Let $\rho$ be as defined in Lemma \ref{denseH}. Let $v \in C_0^{\infty} (\Omega)$ and $v_{\varepsilon}(x) = \rho (|x|/\varepsilon)\, v(x)$. Then
\begin{eqnarray} \label{H1v}
||v_{\varepsilon} - v||^2_{\widetilde {\cal H}} &=& \int_{A_\varepsilon} |x|^{-(N-2)}\, \left ( \nabla_x \rho \left ( \frac{|x|}{\varepsilon} \right ) \right )^2\, v^2\, dx \nonumber \\ && + \int_{A_\varepsilon} |x|^{-(N-2)}\, (\rho-1)\, \nabla_x \rho \left ( \frac{|x|}{\varepsilon} \right ) \cdot \nabla v^2\, dx \nonumber \\ && + \int_{B_{2\varepsilon}} |x|^{-(N-2)}\, (\rho-1)^2\, |\nabla v|^2\, dx,
\end{eqnarray}
where $A_{\varepsilon} = \{x \in \Omega, \varepsilon < |x| < 2\varepsilon  \}$. Integrating by parts the second integral and noting that $\rho'(2) =0$ and $\rho'(1) =0$, we have that
\begin{eqnarray*}
\int_{A_\varepsilon} |x|^{-(N-2)}\, (\rho-1)\, \nabla_x \rho \left ( \frac{|x|}{\varepsilon} \right ) \cdot \nabla v^2\, dx &=& - \int_{A_\varepsilon} |x|^{-(N-2)}\, \left ( \nabla_x \rho \left ( \frac{|x|}{\varepsilon} \right ) \right )^2\, v^2\, dx \\ &&- \varepsilon^{-2} \int_{A_\varepsilon} |x|^{-(N-2)}\, (\rho-1)\, \rho''\, v^2\, dx \\ &&- \varepsilon^{-1} \int_{A_\varepsilon} |x|^{-(N-1)}\, (\rho-1)\, \rho'\, v^2\, dx.
\end{eqnarray*}
Thus, we obtain that
\begin{eqnarray} \label{H1.4}
||v_{\varepsilon} - v||^2_{\widetilde {\cal H}} &=& - \varepsilon^{-1} \int_{A_\varepsilon} |x|^{-(N-1)}\, \left ( \frac{|x|}{\varepsilon} \rho'' + \rho' \right )\, (\rho-1)\, v^2\, dx \nonumber \\ && + \int_{A_\varepsilon} |x|^{-(N-2)}\, (\rho-1)^2\, |\nabla v|^2\, dx.
\end{eqnarray}
Taking $\varepsilon < < 1$, we have that $v^2(x) \sim v^2 (0)$, $x \in A_{\varepsilon}$. Then, for $t=r/\varepsilon$, the first integral in (\ref{H1.4}) is equal to
\[
- \varepsilon^{-1} \int_{A_\varepsilon} |x|^{-(N-1)}\, \left ( \frac{|x|}{\varepsilon} \rho'' + \rho' \right )\, (\rho-1)\, v^2\, dx = c\, \int_{1}^{2} \left ( t\, \rho'(t) \right )'\, \left ( \rho(t)-1 \right )\, dt + O(\varepsilon).
\]
Integrating by parts the last integral and using again that $\rho'(2)=\rho'(1) =0$, we obtain that it is equal to
\begin{equation} \label{H1.5}
- \varepsilon^{-1} \int_{A_\varepsilon} |x|^{-(N-1)}\, \left ( \frac{|x|}{\varepsilon} \rho'' + \rho' \right )\, (\rho-1)\, v^2\, dx = c\, \int_{1}^{2} t\, \left ( \rho'(t) \right )^2\, dt + O(\varepsilon).
\end{equation}
On the other, for $\varepsilon \downarrow 0$, the last integral of (\ref{H1.4})
\[
\int_{B_{2\varepsilon}} |x|^{-(N-2)}\, (\rho-1)^2\, |\nabla v|^2\, dx \leq c\, \int_0^{2\varepsilon} r\, dr \to 0.
\]
Hence, from (\ref{H1.4}) and (\ref{H1.5}) we conclude that
\begin{eqnarray} \label{H1.7}
\lim_{\varepsilon \downarrow 0} ||v_{\varepsilon} - v||^2_{\widetilde {\cal H}} &=& c\, \int_{1}^{2} t\, \left ( \rho'(t) \right )^2\, dt,
\end{eqnarray}
which is a positive quantity. The proof is thus completed. $\blacksquare$ \vspace{0.2cm}

The special cutoff functions that are dense in $\widetilde {\cal H}$ are the ones that allow to prove that $\{0\}$ has zero capacity in two space dimensions. The connection of transformation (\ref{trans}) with the reduction of dimension from $N$ to $2$ is one of the leitmotifs of paper \cite{bv97}. More precisely, we prove the following:
\begin{lemma} \label{densetildeH}
The $C_0^{\infty} (\Omega \backslash \{0\})$-functions are dense in $\widetilde {\cal H}$.
\end{lemma}
\emph{Proof} By the definition of $\widetilde {\cal H}$, it
suffices to prove that we can approximate smooth
functions vanishing at $\partial \Omega$ by smooth functions that
vanish also near $x=0$:
\[
C_0^{\infty} (\Omega) \subset \overline{C_0^{\infty} (\Omega \backslash \{0\})}^{||.||_{\widetilde {\cal H}}}
\]
Let $\rho_\ve(t)$, $\ve\in(0,1)$, be the  cutoff
functions defined as follows:  (i) $\rho_\ve(t)=1$ for $t \geq
\ve$; (ii) $\rho_\ve(t)=0$ for $0<t\leq \ve^2$; (iii) in the
remaining region  $1/\ve^2\le t\le 1/\ve$, it has the special form
$$
\rho_\ve(t)=c_\ve \log(t/\ve^2), \qquad c_\ve=(\log (1/\ve))^{-1}.
$$
For a fixed $v \in C_0^{\infty} (\Omega)$ we define
$v_{\varepsilon}(x) = \rho_\ve (|x|)\, v(x)$.  Note that
$$
|\nabla \rho_\ve(|x|)|=\frac{c_\ve}{|x|} \qquad \mbox{ for } \
x\in A_{\varepsilon} = \{x \in \Omega: \ \varepsilon^2 < |x| <
\varepsilon  \}
$$
being zero otherwise.  Then we have that
\begin{eqnarray}\label{H1a}
||v_{\varepsilon} - v||^2_{W_0^{1,2} (d\mu)} &\leq&
2\int_{A_\varepsilon} |x|^{-(N-2)}\, |\nabla \rho_\ve (|x|)|^2\,
v^2\, dx \nonumber \\ && + 2 \int_{B_{\varepsilon}} |x|^{-(N-2)}\,
(1-\rho_\ve)^2\, |\nabla v|^2\, dx,
\end{eqnarray}
We will prove that letting $\varepsilon \downarrow
0$, the integrals in (\ref{H1a}) go to zero. The only one that is
delicate is the first. We have
\begin{equation} \label{h2}
\dint_{A_\varepsilon} |x|^{-(N-2)}\, |\nabla \rho_\ve (|x|)|^2\, v^2\, dx \le C\|v\|_\infty^2 \int_{\ve^2}^{\ve}c_\ve^2\frac{dr}{r}= C\|v\|_\infty^2(\log (1/\ve))^{-1}
\end{equation}
and this tends to zero as $\ve\to 0$. Approximating the functions $v_{\varepsilon}$, which vanish around the origin, by $C_0^{\infty} (\Omega \backslash \{0\})$  is now standard. $\blacksquare$ \vspace{0.2cm} \\
From Lemmas \ref{denseH} and \ref{densetildeH} we have that the spaces $H$ and $\H$ may both be defined as the closure of $C_0^{\infty} (\Omega \backslash \{0\})$-functions, with respect to (\ref{oldH}) and (\ref{newH}), respectively. However, for such functions these two norm are equal, i.\,e.,
\begin{equation} \label{H1=I}
||w||^2_{\H} = I_{\Omega}(w),
\end{equation}
for any $w \in C_0^{\infty}(\Omega\backslash \{0\})$. Thus,
\begin{proposition} \label{propH=H}
The spaces $H$ and $\H$ coincide.
\end{proposition}
For the space $\widetilde {\cal H}$, which is defined by (\ref{weigted space}), we have that
\begin{lemma} \label{lemmacalHcontains}
The space $\widetilde {\cal H}$ contains all the functions that satisfy $v|_{\partial \Omega} =0$
and $||v||_{\widetilde {\cal H}} < \infty$.
\end{lemma}
\emph{Proof} We first note that the bounded functions with $v_{\partial \Omega} =0$
and $||v||_{\widetilde {\cal H}} < \infty$ belong in $\widetilde {\cal H}$. This follows using the $v_{\varepsilon}$ functions defined in Lemma \ref{densetildeH}, see (\ref{h2}).

Next, for any $v$ with $v_{\partial \Omega} =0$ and $||v||_{\widetilde {\cal H}} < \infty$, we will prove that $v$ is approximated by a sequence of bounded functions in $\widetilde {\cal H}$. Indeed, if we define $v_n$ as follows
\[
v_n(x)=v(x),\;\; \mbox{if}\;\; |v(x)| \le n,\;\;\; v_n(x)=n,\;\; \mbox{if}\;\; |v(x)|> n,
\]
we have that
\begin{equation}
\int_{\Omega} |x|^{-(N-2)}\, |\nabla (v-v_n)|^2\, dx = \int_{C_n}
|x|^{-(N-2)}\, |\nabla v|^2\,dx < \infty.
\end{equation}
Now, it is clear that the sets $C_n=\{x\in \Omega: \ |v(x)|>n \}$
form a monotone family and the  measure tends to zero in the limit
$n\to \infty$. This means that the above integral goes to zero as
$n\to\infty$ and the proof is complete. $\blacksquare$
\vspace{0.2cm}

We have also obtained the above results  following a different
approach. This is included in the Appendix since some readers may
be interested in that approach.

\subsection{On the norm of $H$} \label{sub.normH}

Let us examine further the definition of the norm $N$ that will be considered for the space $H=\cal H$.
We know that the norms $I_\Omega^{1/2}$ and $N$ coincide on functions $H^1_0 (\Omega)$, and also that $I_\Omega(u)$
is larger than $N^2(u)=\|u\|_H^2$ when they differ. More precisely,
 \begin{equation} \label{new norm H}
||u||^2_H = \lim_{\varepsilon \to 0} \left ( I_{B^c_\varepsilon} [u] - \L_{\varepsilon} (u) \right ).
\end{equation}
Now, if for any $u\in H$ we consider a sequence of cutoff approximations $u_\ve(x)=\rho_\ve(x)u(x)$ with $\rho_\ve$ as in Lemma (\ref{densetildeH}), then $\|u_\ve\|_H=I_\Omega(u_\ve)^{1/2}$ and $u_\ve\to u$ in $H$. The limit value
$\|u\|_H^2 =\lim_{\ve\to 0} \|u_\ve \|_H^2$ is what we call the {\sl cutoff value of the Hardy functional}, and produces the correct norm on $H$.

As a conclusion, the space $H$ as it is defined in Vazquez-Zuazua
\cite{vz00} is a well defined space, as the completion of
$C_0^{\infty} (\Omega)$-functions with respect to the norm
$||\phi||^2_H = I_{\Omega}(\phi)$, $\phi \in C_0^{\infty}
(\Omega)$. In this space there exist ``bad functions'', such that
$I_{\Omega}(u)$ defined as an improper integral does not coincide
with the limit of the sequence of cutoff approximations. Even
more, it can happen that the principal value of the integrals in
$I_\Omega$ is not well defined, either oscillating or infinite.
For example, let $u$ behaving at the origin as $|x|^{-(N-2)/2}$;
for this function, the quantity $I_{\Omega}(u)$ is well defined,
but its norm in $H$ is not $I_{\Omega}(u)$, but it is equal to
$I_{\Omega}(u)-\L(u)$. For the first eigenfunction $e_1$: the norm
is not $I_{\Omega}^{1/2}(e_1)$, but
\[
||e_1||^2_H = I_{\Omega}(e_1) - \L(e_1) = I_{\Omega}(e_1) - \frac{N(N-2)}{2}\, \omega_N = \mu_1.
\]

As a result the minimization problem (\ref{01}) and the following one
\begin{equation}\label{01a}
\min_{u \in H} \frac{I_{\Omega}(u)}{||u||^2_{L^2}}
\end{equation}
are not the same, since the first admits a minimizer while the other does not. This will also happen in the case of the minimizers of (\ref{otherminimizers}), see Subsection \ref{sub.furtherprop}, and in the case of the minimizer of the Improved Hardy-Sobolev inequality in the radial case, \cite{z10}.

As it is stated in \cite[pg. 127]{vz00}, the solutions of the
stationary and evolution problem in the critical case, $c_* =
(N-2)^2/4$, can be obtained as the limit when $c \uparrow c_*$ of
the solutions of the subcritical case $c<c_*$. This is true in the
space $H$, with norm (\ref{new norm H}), and not through $I$. For
example, let $e_{1,c}$ be the eigenfunctions corresponding to the
principal eigenvalues of the subcritical case, on the unit ball.
In \cite[Theorem 3.3]{vz00} the exact formula of these
eigenfunctions were obtained:
\[
e_{1,c} (r) = r^{-(N-2)/2}\, J_m (z_{m,1}\, r),
\]
where $m^2 = c_* - c >0$ and $z_{m,1}$ is the first zero of the
Bessel function $J_m$. Then
\[
||e_{1,c}||^2_{H} = I_{\Omega}(e_{1,c} (r)) \to I_{\Omega}(e_1) - \L(e_1) = ||e_1||^2_H,
\]
as $c \uparrow c_*$.

This strange situation has to do with the definition of
$I_{\Omega}(u)$; the existence of the difference of two improper
(and comparable) integrals. This is the difference between the
case of a classical Sobolev space or the case of the Hardy
functional with constant less than the best constant.

%
%
\section{Application to the evolution problem} \label{sec.evol}
\setcounter{equation}{0}
We now justify that  the results described in  \cite{vz00} for the
solutions of Problem (\ref{eq1}) hold in the space $\H$ defined in
(\ref{newH}), with $\Omega$ a bounded domain of $\mathbb{R}^N$, $N
\geq 3$. Is clear that this space is a Hilbert space and, as
stated in \cite[Theorem 2.2]{vz00},  $\H$ is imbedded continuously
in the Sobolev space $W^{1,q}(\Omega)$, $1 \leq q <2$, i.e.,
\begin{equation}\label{embedding}
||u||_{\H} \geq C(q,\Omega)\, ||u||_{W^{1,q}(\Omega)},\;\;\; 1 \leq q <2.
\end{equation}
Thus, the compact imbedding
\begin{equation}\label{compLp}
\H \hookrightarrow L^p (\Omega),\;\;\; 1\leq p < \frac{2N}{N-2}.
\end{equation}
holds. Moreover, we may justify that all the results concerning
the spectrum of the related eigenvalue problem given in
\cite{vz00} hold for $\H$. For this it is sufficient to proceed as
follows. Assume the Cauchy-Dirichlet problem for equation
\begin{equation} \label{Cauchy-Dirichlet_u}
u_t = \Delta u + \biggr ( \frac{N-2}{2} \biggr
)^2 \frac{u}{|x|^2}.
\end{equation}
and the corresponding eigenvalue problem
\begin{equation} \label{eigen_u}
\Delta u + \biggr ( \frac{N-2}{2} \biggr
)^2 \frac{u}{|x|^2} + \mu\, u=0.
\end{equation}
Setting $u = |x|^{-\frac{N-2}{2}}v$, (\ref{Cauchy-Dirichlet_u})
becomes equation (\ref{Cauchy-Dirichlet_v}) away from zero. The
distributional operator $L$ given by
\begin{equation} \label{L}
L(v) = - |x|^{N-2}\; \nabla \cdot \left( |x|^{-(N-2)}\, \nabla v \right)
\end{equation}
is an isomorphism from $\widetilde {\cal H}$ into its dual
$\widetilde {\cal H}^{-1}$ and from (\ref{compLp}) we have the
triplet
\[
\widetilde {\cal H} \to \widetilde X \to \widetilde {\cal H}^{-1}.
\]
Then, by restriction, we may define the surjective operator $L_* : D(L_*) \subset \widetilde X \to \widetilde X$, with domain
\[
D_* = \left \{ f \in \widetilde {\cal H}: |x|^{N-2}\; \nabla \cdot
\left( |x|^{-(N-2)}\, \nabla f \right) \in \widetilde X \right \}.
\]
If we denote $L_*$ with $L$, is easy to see that $L$ is self-adjoint with compact inverse. So, we can form
an orthonormal basis of eigenfunctions in $\widetilde X$, with eigenvalue sequence
\begin{equation} \label{eigenv}
0 < \mu_1 \leq \mu_2 \leq ... \leq \mu_k \leq ... \to \infty.
\end{equation}
(As in \cite{vz00}, separation of variables allows to find explicit formula for the eigenfunctions and eigenvalues when $\Omega$ is a ball). Returning now to Problem (\ref{Cauchy-Dirichlet_u}), we have that the result of \cite{vz00}, concerning the spectrum of the associated eigenvalue problem as well as the asymptotic behavior of the evolution problem, both are true.

\subsection{Energy calculations} \label{sub.energycalc}

In order to justify the above choice we make the following
calculations. Let $\phi_k$ be the eigenfunctions corresponding to
(\ref{eigenv}). Assume that $v(t,x) $ is a solution of
(\ref{Cauchy-Dirichlet_v}) and $v \in L^2(0,T:\widetilde {\cal
H})$. If $v(t,x) = c_k (t)\, \phi_k$, we have that
\[
||v(t)||^2_{\widetilde X} = \sum_{k=1}^{\infty} c_k^2 (t)\;\;\;
\mbox{and}\;\;\; ||v(t)||^2_{\widetilde {\cal H}} =
\sum_{k=1}^{\infty} c_k^2 (t)\, ||\phi_k||^2_{\widetilde {\cal
H}}.
\]
If we denote by $E (t) := ||v(t)||^2_{\widetilde X}$, we obtain
that
\begin{equation}\label{2.2.1}
\frac{dE(t)}{dt} = 2\, \sum_{k=1}^{\infty} c_k (t)\, c'_k (t).
\end{equation}
and for $E_{v \phi_k} (t) := <v(t),\phi_k>_{\widetilde X}$ we get that
\begin{equation}\label{2.2.2}
\frac{dE_{v \phi_k} (t)}{dt} = c'_k (t).
\end{equation}
However, the evolution equation implies that
\begin{equation}\label{2.2.3}
\frac{dE_{v \phi_k} (t)}{dt} = - <v,\phi_k>_{\widetilde {\cal H}}=
-c_k\, ||\phi_k||^2_{\widetilde {\cal H}}.
\end{equation}
Thus, (\ref{2.2.2}) and (\ref{2.2.3}) imply that $c'_k (t) = -c_k
||\phi_k||^2_{\widetilde {\cal H}}$. Then, from (\ref{2.2.1}) we
conclude that
\begin{equation}\label{2.2.4}
\frac{dE(t)}{dt} = - 2\, \sum_{k=1}^{\infty} c^2_k (t)\,
||\phi_k||^2_{\widetilde {\cal H}} = -2 ||v||^2_{\widetilde {\cal
H}} = -2\, \int_{\Omega} |x|^{-(N-2)}\, |\nabla v|^2\, dx.
\end{equation}

$\bullet$ Let us now see where the original Hardy functional fails
by performing this energy calculation in terms of $u$. To avoid
the singularity we calculate the energy in
$\Omega_\ve=\Omega\setminus B_\ve(0)$with $\ve>0$ small.
$$
\begin{array}{l}
\dfrac{d}{dt}\dint_{\Omega_\ve}u^2\,dx=\dint_{\Omega_\ve}u\,\Delta
u\,dx+c_* \dint_{\Omega_\ve}\dfrac{u^2}{|x|^2}\,dx\\
=- \dint_{\Omega_\ve}\left(|\nabla
u|^2-c_*\dfrac{u^2}{|x|^2}\right)\,dx + \int_{|x|=\ve} u u_n \,dS.
\end{array}
$$
Using $u=v\,r^{-(n-2)/2}$n with $r=|x|>0$, we calculate the last
integral
$$
\int uu_n\,dS=-\dfrac{n-2}{2}\ve^{1-n}\dint_{r=\ve} v^2dS+
\ve^{2-n}\dint_{r=1} v v_r dS.
$$
Recalling now formula (\ref{twoint}) we arrive at
$$
\dfrac{d}{dt}\dint_{\Omega_\ve}u^2\,dx=-||v||^2_{\widetilde {\cal
H} (B^c_\varepsilon)} + \ve^{2-n}\dint_{r=\ve} vv_r dS
$$
If we now pass to the limit $\ve\to0$ and realize that the last
integral goes to zero whenever $v\in \widetilde{\cal H}$, we get
formally
$$
\dfrac{d}{dt}\dint_{\Omega}u^2\,dx=-||v||^2_{\widetilde {\cal H}
(B^c)}
$$
This is precisely the choice made by our setting. Let us briefly
justify the term that went to zero in average when $v$ is bounded:
integrating in $\ve$ from $\delta$ to $2\delta$ and putting
$A_\delta=\{ x: \delta\le |x| \le 2\delta \}$, we have the
following estimate for the average
$$
|\frac1{\delta}\int_{\delta}^{2\delta} \ve^{2-n}d\ve \dint_{r=\ve}
vv_r dS| \le \frac{C}{\delta} \int_{A_\delta} |\nabla
v|\,|x|^{2-n}dx\le \frac{C}{\delta} \left(\int_{A_\delta} |\nabla
v|^2\,|x|^{2-n}dx\right)^{1/2}\left(\int_\delta^{2\delta}rdr\right)^{1/2}
$$
that goes to zero as $\delta\to 0$. So, the weak formulation (or
the energy equation) of (\ref{eq1}) is the following
\[
\frac{1}{2} \int u_t^2 = -||u||^2_H = - \lim_{\varepsilon \to 0} \left ( I_{B^c_\varepsilon} [u] - \L_{\varepsilon} (u) \right )
\]
for every $u \in H$. The space $H$ is really the energetic space.

\medskip

\section{Further properties of the spaces} \label{sec.h}
\setcounter{equation}{0}

\subsection {Further properties of  $\H$ and $\widetilde {\cal H}$} \label{sub.furtherprop}

We observe that imbedding (\ref{embedding}) gives the following
corollary, which completes the results obtained in \cite{ckn} (see
also \cite{cw01}) concerning the
Caffarelli--Kohn--Nirenberg Inequalities, in a bounded domain, in
the limiting case where $a=\frac{N-2}{2}$.
\begin{corollary} {\rm (Critical Caffarelli--Kohn--Nirenberg Inequalities)}
Assume that $v_n$ is a bounded sequence in $\widetilde {\cal H}$.
Then $u_n = |x|^{-(N-2)/2}\, v_n$ is a bounded sequence in $\H$.
The compact imbeddings {\rm (\ref{compLp})} imply that, up to some
subsequence, $u_n$ converges in $L^{q}(\Omega)$ to some $u$. Thus,
we obtain the compact imbeddings
\begin{equation}\label{compactimbed1}
\widetilde {\cal H} \hookrightarrow L^q
(|x|^{-q(N-2)/2}dx,\Omega),\;\;\; \mbox{for any}\;\;\; 1 \leq q <
\frac{2N}{N-2}.
\end{equation}
Then, for every $0 \leq s \leq \frac{N-2}{2}q$, we further obtain the compact imbeddings
\begin{equation}\label{compactimbed2}
\widetilde {\cal H} \hookrightarrow L^q
((|x|^{-s}dx,\Omega),\;\;\; \mbox{for any}\;\;\; 1 \leq q <
\frac{2N}{N-2},
\end{equation}
where the weighted space $L^q (w(x)\,dx,\Omega)$ is defined as the
closure of $C_{0}^{\infty}(\Omega)$ functions under the norm
\[
||\phi||_{L^q (w(x),\Omega)} = \left ( \int_{\Omega} w(x)\,
|\phi|^q\, dx \right )^{\frac{1}{q}}.
\]
\end{corollary}
\begin{remark}
In {\rm (\ref{compactimbed1})} it is clear that $q$ cannot reach
$\frac{2N}{N-2}$. For this value of $q$, the best that we can have
are Improved Hardy-Sobolev Inequalities, see \cite{z10} and the
references therein.
\end{remark}

In addition, we can relate these spaces, in the radial case, with
the space $\D^{1,2}(\mathbb{R}^N)$, which is defined as the closure
of $C_{0}^{\infty}(\mathbb{R}^N)$ functions under the norm
\[
||\phi||^2_{D^{1,2}(\mathbb{R}^N)} = \int_{\mathbb{R}^N} |\nabla
\phi|^2\, dx.
\]
If we denote by $\H_r (\Omega)$, $\widetilde{\cal H}_{r}$ and
$\D_r^{1,2}(\mathbb{R}^N)$ the subspaces of $\H $, $\widetilde
{\cal H}$ and $\D^{1,2}(\mathbb{R}^N)$, respectively, which
consist of radial functions, we have that
\begin{proposition} \label{W=D}
For some function $v \in \widetilde{\cal H}_{r}(B_R)$ we set
\begin{equation} \label{transxw}
v(|x|) = w(t),\;\;\;\;\;\; t = \left ( - \log \left (
\frac{|x|}{R} \right ) \right )^{-\frac{1}{N-2}}.
\end{equation}
Then, $v \in \widetilde{\cal H}_{r}(B_R)$ if and
only if $w \in \D_r^{1,2}(\mathbb{R}^N)$ and
\begin{equation}\label{eqWD}
||v||_{\widetilde{\cal H}_{r}(B_R)} = (N-2)^{-1}\, ||w||_{D_r^{1,2}(\mathbb{R}^N)}.
\end{equation}
\end{proposition}
Observe that (\ref{eqWD}) is independent of the radius $R$ and in the case where $N=3$ the norm in $\widetilde{\cal H}_{r}(B_R)$ coincides with the norm in $\D_r^{1,2}(\mathbb{R}^N)$. Moreover, the definition of $\H$ and (\ref{eqWD}) imply that
\begin{corollary}
For some function $u$ we set
\begin{equation} \label{transx}
w(t) = |x|^{\frac{N-2}{2}} u(|x|), \;\;\;\;\;\; t = \left ( - \log
\left ( \frac{|x|}{R} \right ) \right )^{-\frac{1}{N-2}}.
\end{equation}
Then, $u \in \H_r (B_R)$ if and only if $w \in \D_r^{1,2}(\mathbb{R}^N)$ and
\begin{equation}\label{eqHD}
||u||^2_{\H_r (B_R)} = (N-2)^{-1}\, ||w||^2_{\D_r^{1,2}(\mathbb{R}^N)}.
\end{equation}
\end{corollary}

Transformation (\ref{transxw}) was used in \cite{z10}. For a discussion concerning the construction of such transformations for general Hardy inequalities we refer to \cite{vzog}.

\subsection {Nonexistence of $H_0^1$-minimizers}
Transformation
(\ref{transx}) provides us with an extra argument concerning the
nonexistence of $H_0^1$-minimizers. And not only this; we are able
to obtain the exact behavior of these minimizers at the
singularity. We will prove that these minimizers belong to $H$,
they do not belong to $H_0^1$ and their behavior at the origin is
exactly $|x|^{-(N-2)/2}$.

As an example, we will prove that the minimizer of (\ref{01}) cannot be an
$H_0^1$-function. As already mentioned in the introduction, this
minimizer is the function $e_1$, given by (\ref{e1}), and does not
belong to $H_0^1$.

Assume to the contrary that $u \in H_0^1$ is a minimizer of
(\ref{01}). Then, $u$ may be chosen to be a nonnegative and
radial function, i.e. $u(x)=u(r) \geq 0$. Let $w$ be the
transformation of $u$, through (\ref{transx}). Since $u \in H_0^1$
we obtain that
\begin{equation}\label{nonex1}
w(0)=0.
\end{equation}
Moreover, we have that $w \in D^{1,2}(\mathbb{R}^N)$ is a minimizer of
\begin{equation}\label{nonex2}
\frac{1}{(N-2)^2}\, \frac{\int_{\mathbb{R}^N} |\nabla w|^2\, dx}{\int_{\mathbb{R}^N} V(|x|)\, w^2\,
 dx}\,,
\end{equation}
where $V(|x|) = |x|^{-2(N-1)}\, e^{-2\, |x|^{-(N-2)}}$. Note that if we set $V(0)=0$, $V$ is a continuous function.
Then, $w$ should be a nonnegative solution of the Euler-Lagrange equation corresponding to (\ref{nonex2}):
\[
- \Delta w = c(N)\, V(|x|)\, w,\;\;\; w \in D^{1,2}(\mathbb{R}^N).
\]
However, the application of the Maximum Principle contradicts (\ref{nonex1}), hence (\ref{01}) cannot have a $H_0^1$-minimizer.

This argument may be applied to more general problems;
\begin{proposition} \label{propNonexist}
Let $\Omega$ be a bounded domain of $\mathbb{R}^N$, $N \geq 3$,
containing the origin. Then, minimizers of
\begin{equation}\label{minproblems}
\min_{u \in H} \frac{||u||^2_{H(\Omega)}}{\left ( \int_{\Omega} |u|^{q}\, dx \right )^{2/q}},\;\;\;\;\; 1 \leq q < \frac{2N}{N-2},
\end{equation}
cannot exist in $H_0^1 (\Omega)$.
\end{proposition}
\emph{Proof} The compact imbeddings (\ref{compactimbed2}) imply that the minimization problems (\ref{minproblems})
have a solution, let us denote it by $c_{\Omega, q}$, and this solution is attained by some function $u_{\Omega, q}$.

If $\Omega$ is radially symmetric, then the conclusion follows
using the same argument that was applied for $e_1$. In the case
where $\Omega$ is an arbitrary domain, we make a symmetrization
that replaces $\Omega$ by a ball $B_R$ with the same volume and
the function $u_{\Omega, q}$ by its symmetric rearrangement
$u^*_{\Omega, q}$.

Assume that $u_{\Omega, q} \in H_0^1 (\Omega)$, then $||u_{\Omega,
q}||^2_{H(\Omega)} = I_{\Omega} [u_{\Omega, q}]$. It is well known
that the rearrangement does not change the $L^q$-norm, increases
the integral $\int u^2/|x|^2\, dx$. Also, decreases the $H_0^1
(\Omega)$- norm, so that $u^*_{\Omega, q}$ belongs to $H_0^1
(B_R)$. Finally, we get the contradiction that $u^*_{\Omega, q}$
is the minimizer $u_{B_R, q}$ and belongs to $H_0^1 (B_R)$. Thus,
$u_{\Omega, q}$ cannot exist in $H_0^1 (\Omega)$. $\blacksquare$
\vspace{0.2cm}

The case $q=\frac{2N}{N-2}$ has the same quantitative behavior (in the radial case) and this maybe obtained following the same argument.

As mentioned before, the case of $e_1$ and the minimizer of the
improved Hardy-Sobolev Inequality (in the radial case), belong to
the second case of Subsection \ref{sec.hf}. This means that they
behave at the origin like $|x|^{-(N-2)/2}$. Thus, the Hardy
functional for these functions is a well defined positive number,
although it does not represent their $H$-norm. These functions do
not belong to the "worst" cases, where $I_{\Omega}$ is not well
defined or is infinite. As a corollary of the previous argument we
have that the same happens to every minimizer $u_{\Omega, q}$ of
(\ref{minproblems}).
\begin{corollary} \label{cornonex}
Every minimizer $u_{\Omega, q}$ of {\rm (\ref{minproblems})}
behaves at the origin like $|x|^{-(N-2)/2}$.
\end{corollary}
\emph{Proof} Assume the opposite. Since $u_{\Omega, q}$ does not
belong in $H_0^1$, in the sense of Subsection \ref{sec.hf}, we
have that the Hardy functional of $u_{\Omega, q}$ is not well
defined or is infinite:
\[
I_{B^c_\varepsilon} [u_{\Omega, q}]\;\;\;\;\;\; \mbox{oscillates or is infinite, as}\;\; \varepsilon \to 0.
\]
Since, $u_{\Omega, q} \in H(\Omega)$, from (\ref{new norm H}) we have that the same happens to $\L_{\varepsilon} (u_{\Omega, q})$, as $\varepsilon \to 0$.

We consider first the radial case, i.e. $u_{\Omega, q}$ is a
radial function and $\Omega = B_R$. Let $w_{\Omega, q}$ be the
transformation of $u_{\Omega, q}$, through (\ref{transx}). Then,
from the definition (\ref{Lepsilon}) of $\L_{\varepsilon}$,
$w_{\Omega, q}(0)$ is not well defined or is infinite,
respectively. On the other, $w_{\Omega, q} \in
D^{1,2}(\mathbb{R}^N)$ is a minimizer of
\[
\frac{1}{(N-2)^{1+1/q}}\, \frac{\int_{\mathbb{R}^N} |\nabla w|^2\, dx}{\left ( \int_{\mathbb{R}^N} V_q(|x|)\, w^q\, dx \right )^{2/q}},
\]
i.e., $w_{\Omega, q}$ is the solution of the corresponding Euler-Lagrange equation
\[
- \Delta w = c(q,N)\, V_q(|x|)\, |w|^{q-2}\, w,\;\;\; w \in D^{1,2}(\mathbb{R}^N),
\]
with $V_q (|x|)= |x|^{-2(N-1)}\, e^{-f(q)\, |x|^{-(N-2)}}$, $f(q) = N-(N-2)q/2 >0$. $V_q$ is a continuous function, if we set $V_q(0)=0$. However, standard regularity theory implies that this is a contradiction.

Next we treat the general case. We consider $u^*_{\Omega, q}$, the symmetric rearrangement of $u_{\Omega, q}$.
In the case where $\lim_{x \to 0} |x|^{(N-2)/2}\, u_{\Omega, q} (x) = \infty$, $u^*_{\Omega, q}$ does the same at zero, since their level sets have the same volume. In the case, where $\lim_{x \to 0} |x|^{(N-2)/2}\, u_{\Omega, q} (x)$ oscillates, \color{blue} we make the \emph{unproved hypothesis} that $u^*_{\Omega, q}$ cannot be defined at $0$. \normalcolor
However, $u^*_{\Omega, q}$ should be a radial minimizer and this is a contradiction. $\blacksquare$ \vspace{0.2cm} \\
As in the case of $e_1$, we emphasize the fact that the minimization problem (\ref{minproblems}) and the following
\[
\min_{u \in H} \frac{I_{\Omega}(u)}{||u||^{2/q}_{L^q}}
\]
are not the same.

%

\section{The case of the exterior domain} \label{sec.exd}
\setcounter{equation}{0}

We consider Problem (\ref{exd1})  describing the evolution (up to some weight) of the Hardy potential in an exterior domain. We may fix $\delta=1$. Our arguments will be based on the unitary equivalence with the previous problem posed on a ball. For that we use the Kelvin transform in the form
\begin{equation}\label{kelvintrans}
u(x) = |y|^{N-2}\, w(y),\;\;\;\;\;\; y=\frac{x}{|x|^2}\,.
\end{equation}
These formulas  transform solutions $u(x,t)$ of Problem (\ref{eq1}), i.\,e.,
\begin{equation}\label{exd2}
\left\{\begin{array}{ccc}
u_t (x,t) &=& \Delta u (x,t)+ c_*\, \dfrac{u(x,t)}{|x|^2},\; x \in B_1 (0),\; t>0,\\
u(x,0) &=& u_0(x),\;\; \mbox{for}\;\; x \in B_1 (0), \nonumber \\ u(x,t)&=&0\;\;\mbox{in}\;\;\partial B_1 (0),\;t>0\,.
\end{array}\right.
\end{equation}
defined in the unit sphere $B_1 (0)$ into solutions $w(y,t)$ of Problem (\ref{exd1}) posed in $B^c_1 (0)$, its complement in $\mathbb{R}^N$.
We will write the transformation as $u=\K(w)$. Note that for smooth functions holds that
\[
\Delta_x u(x) = |y|^{N+2}\, \Delta_y w(y)\;\;\;\;\;\; \mbox{and}\;\;\;\;\;\; \frac{u(x)}{|x|^2} = |y|^{N+2}\, \frac{w(y)}{|y|^2}.
\]
The equivalence of the equations is clear equivalence for $x \ne 0$. The differences will appear in the energy near the singularity versus the ``energy at infinity''.

\medskip

\noindent {\bf Basic properties.} We will address the questions of existence of solutions for Problem (\ref{exd1}) by means of this equivalence, that will also be used as a clue to the proper variational formulation that will be followed.
As we saw, the functional space which corresponds to (\ref{exd2}) is $H$ with norm given in (\ref{new norm H}). Our
proposal is to use this formulation for the definition of the new space, $\W$.
The space $\W$ is defined as the isometric space of $H$ under the Kelvin transformation  (\ref{kelvintrans}). In other words, $\W$ is defined as the
completion of the set
\[
\left \{ w (y) = |y|^{-N+2}\, u \left( \frac{y}{|y|^2} \right ),\;\;\; u \in C_{0}^{\infty}
(B_1 (0)),\;\;\; |y| \geq 1 \right \},
\]
under the norm $\|w\|_{\W}$ defined by
\begin{equation} \label{normWH}
||w||^2_{\W} = \lim_{\varepsilon \to 0} \left ( I_{B_1 (0) \backslash B_\varepsilon} [u] - \L_{\varepsilon} (u) \right ),\;\;\;\;\;\; u=\K(w).
\end{equation}

The first eigenpair of the corresponding eigenvalue problem
\begin{equation}\label{exdeigenproblem}
- \Delta w - c_* \frac{w}{|y|^2} =\mu\, |y|^{-4}\, w
\end{equation}
is
\begin{equation}\label{exdeigenpair}
\mu_1 = z^2_{0,1},\;\;\;\;\; \tilde{e}_{1} = |y|^{-(N-2)/2}\, J_0 \left ( \frac{z_{0,1}}{|y|} \right ).
\end{equation}
We  have the following {\bf Result:} {\sl
The well posedness of {\rm (\ref{exd1})} in the space $\W$ is understood through the unitary equivalence with $H$. The existence, uniqueness and stabilization results of Problem {\rm (\ref{exd2})} apply for Problem {\rm  (\ref{exd1})}.}

\medskip

\noindent {\bf Hardy functional.} Next, we investigate the connection of the space $\W$ with the Hardy functional, $I_{B^c_1 (0)}$ defined as:
\begin{equation} \label{Hardyfunctionalex}
I_{B^c_1 (0)} [\phi] = \int_{\mathbb{R}^N \backslash B_{1} (0)} |\nabla \phi|^2\, dx -
\left ( \frac{N-2}{2} \right )^2 \int_{\mathbb{R}^N \backslash B_{1} (0)}
\frac{\phi^2}{|x|^2}\, dx\,,
\end{equation}
which is positive for any compactly supported $\phi \in C^{\infty} (B^c_1 (0))$ that vanishes on the boundary. We denote by $I_{\varepsilon}$ and by $I_{1/\varepsilon}$ the Hardy functional defined on $B_1 (0) \backslash B_\varepsilon$ and $B_{1/\varepsilon} (0) \backslash B_1 (0)$, respectively.  Recall that  $\L_{\varepsilon}(u)$ is given by {\rm (\ref{Lepsilon})}.

\begin{lemma} We have the following fundamental relation:
\begin{equation}\label{Iwu}
I_{\varepsilon} [u] = I_{1/\varepsilon} [w] + 2 \L_{1/\varepsilon} (w)\,,
\end{equation}
where
\[
\L_{1/\varepsilon} (w) = \frac{N-2}{2}\, \varepsilon\,
\int_{\partial B_{1/\varepsilon} (0)} w^2\, dS\,.
\]
Moreover, it is clear that if $u = \K (w)$ then
\begin{equation}\label{Lwu}
\L_{\varepsilon} (u) = \L_{1/\varepsilon} (w),
\end{equation}
\end{lemma}

\noindent {\sl Proof.} In order to justify Equations (\ref{Iwu}) and (\ref{Lwu})
we proceed as follows. Let $w \in \W$, i.e. there exists $u \in H$, such that $u=\K(w)$. Then,
\begin{eqnarray} \label{exd3}
\int_{B_1 (0) \backslash B_{\varepsilon} (0)} |\nabla_x u|^2\, dx &=& \int_{B_1 (0) \backslash B_{\varepsilon} (0)} |x|^{-2N+4}\, \left |\nabla_x w \left ( \frac{x}{|x|^2} \right ) \right |^2\, dx \nonumber \\
&&+ (-N+2)^2\, \int_{B_1 (0) \backslash B_{\varepsilon} (0)} |x|^{-2N+2}\, w^2 \left ( \frac{x}{|x|^2} \right )\, dx \nonumber \\
&&+ (-N+2)\, \int_{B_1 (0) \backslash B_{\varepsilon} (0)} |x|^{-2N+3}\, \frac{x}{|x|} \cdot \nabla_x w^2 \left ( \frac{x}{|x|^2} \right )\, dx.
\end{eqnarray}
Integrating by parts the last integral we get that
\begin{eqnarray*}
\int_{B_1 (0) \backslash B_{\varepsilon} (0)} |x|^{-2N+3}\, \frac{x}{|x|} \cdot \nabla_x w^2 \left ( \frac{x}{|x|^2} \right )\, dx &=& (N-2)\, \int_{B_1 (0) \backslash B_{\varepsilon} (0)} |x|^{-2N+2}\, w^2 \left ( \frac{x}{|x|^2} \right )\, dx \\
&& - \varepsilon^{-2N+3}\, \int_{\partial B_{\varepsilon} (0)} w^2 \left ( \frac{x}{|x|^2} \right )\, dS.
\end{eqnarray*}
Then from (\ref{exd3}) we have that
\begin{eqnarray} \label{exd4}
\int_{B_1 (0) \backslash B_{\varepsilon} (0)} |\nabla_x u|^2\, dx &=& \int_{B_1 (0) \backslash B_{\varepsilon} (0)} |x|^{-2N+4}\, \left |\nabla_x w \left ( \frac{x}{|x|^2} \right ) \right |^2\, dx \nonumber \\
&&+ (N-2)\, \varepsilon^{-2N+3}\, \int_{\partial B_{\varepsilon} (0)} w^2 \left ( \frac{x}{|x|^2} \right )\, dS.
\end{eqnarray}
However, the determinant of the Jacobian of the Kelvin transformation in $d \geq 2$ dimension is equal to $-|x|^{2d}$ and moreover
\[
\left | \nabla_x w \left ( \frac{x}{|x|^2} \right ) \right |^2 = |x|^{-4}\, |\nabla_y w (y)|^2.
\]
Finally, we conclude that
\begin{eqnarray*}
\int_{B_1 (0) \backslash B_{\varepsilon} (0)} |\nabla_x u|^2\, dx &=& \int_{B_{1/\varepsilon} (0) \backslash B_1 (0)} |\nabla_y w (y)|^2\, dy \\
&& + (N-2)\, \varepsilon\, \int_{\partial B_{1/\varepsilon} (0)} w^2 (y)\, dS.
\end{eqnarray*}
Thus, (\ref{Iwu}) and (\ref{Lwu}) hold. \qed

\medskip

When we apply these results to functions in the class  $\W$ (by density), we are able to give the following unexpected definition of the norm of $\W$,
\begin{equation} \label{normW}
||w||^2_{\W} = \lim_{\varepsilon \to 0} \left ( I_{1/\varepsilon} [w] + \L_{1/\varepsilon} (w) \right ).
\end{equation}
So, the weak formulation (or the energy equation) of (\ref{exd1}) translates into
\[
\frac{1}{2} \int w_t^2 = -||w||^2_{\W} = - \lim_{\varepsilon \to 0} \left ( I_{1/\varepsilon} [w] + \L_{1/\varepsilon} (w) \right ),
\]
for every $w \in \W$.

\medskip

\textbf{Conclusions and remarks.} (i) We have shown that a correcting term  also appears as  in the energy analysis of the problem posed in the exterior domain. Actually, the correcting term has the same absolute value as the Hardy singularity energy
considered in the problem in a bounded domain, but now it represents a kind of energy at infinity. However, there is a big difference from the bounded domain case since in this case the singular energy acts in an additive way to the usual Hardy integral.

\noindent (ii) Moreover, as we will see below, this new term may be the main part of the total energy, since $I_{1/\varepsilon}$ may be also a negative quantity: we do the calculations for $\tilde{e}_1$, given in (\ref{exdeigenpair}),
\[
||\tilde{e}_1||^2_{\W} = ||e_1||^2_{H} = \mu_1 = z^2_{0,1} \simeq 5.76
\]
However,
\[
||\tilde{e}_1||^2_{\W} = I_{\mathbb{R}^N \backslash B_1 (0)} (\tilde{e}_1) + \L(\tilde{e}_1).
\]
Thus,
\[
I_{\mathbb{R}^N \backslash B_1 (0)} (\tilde{e}_1) = 5.76 - \frac{N(N-2)}{2}\, \omega_N,
\]
which is negative for $N=3$.

\noindent (iii) As we saw in Subsection \ref{sec.hf}, $\L_{1/\varepsilon}$ may have a bad behavior; oscillating or tending to infinity. In this cases, the Hardy functional $I_{1/\varepsilon}$ becomes negative and has the same behavior with $\L_{1/\varepsilon}$, so that the sum of them to become a positive real number.

\noindent (v) As a conclusion, we can say that $\L$ is a hidden energy that ``comes" from infinity and is not only a gain of the total energy but it may represent the main part of the total energy. We are not in the position to say that, but we feel that this energy may represent a simple form of ``dark energy".

\noindent (vi)  Assume (\ref{eq1}) posed on whole $\mathbb{R}^N$.  Then, the HS energy is present only at zero and not at infinity, cf.  next section and \cite{vz00} and \cite{vzog}. However,  the Kelvin transformation implies that the corresponding problem for Equation (\ref{exd1}) is well defined on $\mathbb{R}^N$, and it can be checked that the HS energy is present as a dark  energy at infinity, and not at the origin. For the details, we refer to \cite{vzog}.


\section{A Hardy-Poincar\'{e} inequality on $\mathbb{R}^N$} \label{sec.rn}
\setcounter{equation}{0}
Throughout this section we refer to Hardy functional $I_{\mathbb{R}^N}$, defined on $\mathbb{R}^N$, $N \geq 3,$
\begin{equation}\label{hardyfunctionalRN}
I_{\mathbb{R}^N}[u] = \int_{\mathbb{R}^N} |\nabla u|^2\, dx - \biggr ( \frac{N-2}{2} \biggr
)^2 \int_{\mathbb{R}^N} \frac{u^2}{|x|^2}\, dx,
\end{equation}
with limits taken in the sense of principal value if both integrals diverge. First, we discuss the reason why the transformation (\ref{trans}) is not proper in the case of $\mathbb{R}^N$. Let $u \in C_0^{\infty} (\mathbb{R}^N)$ and $v$ given by (\ref{trans}). Then this formula simplifies $I_{\mathbb{R}^N}(u)$ since:
\begin{equation}\label{qf.vrn}
I_{\mathbb{R}^N}(u) = I_{1,\mathbb{R}^N}(v) = \int_{\mathbb{R}^N} |x|^{-(N-2)} |\nabla v|^2\, dx\,.
\end{equation}
However, for radial functions $I_{1,\mathbb{R}^N}(v)$ is equal, up to some constant, to
\[
I_{1,\mathbb{R}^N}(v) = c\, \int_{\mathbb{R}^2} |\nabla v|^2\, dx,
\]
and it is well known that $D^{1,2}(\mathbb{R}^2)$ is a not well defined Hilbert space; it contains the constant functions, so we cannot have an imbedding in $L^q (\mathbb{R}^2)$, for any $1 \leq q$. Thus, defining the space $\widetilde {\cal H} (\mathbb{R}^N)$, as in the bounded domain case through $I_{1,\mathbb{R}^N}(v)$, we get that this is not well defined since it contains the constant functions. This also a way to see that the Hardy functional on $\mathbb{R}^N$ cannot be improved by the $L^2 (\mathbb{R}^N)$ norm, for any $C_0^{\infty} (\mathbb{R}^N)$- function.

\medskip

\textbf{Definition of the new proper transformation} As it will be clear, the transformation (\ref{trans}) is not the only proper one to simplify the Hardy functional. The way we follow to solve the above difficulty is to introduce another convenient variable by means of the formula
\begin{equation} \label{transRN}
u (x)= |x|^{-\frac{N-2}{2}}\, J_0 (|x|)\, v(x):= \T_{\mathbb{R}^N} (v),
\end{equation}
where $J_0$ is the Bessel function. Note that $J_0 (0) =1$, up to normalization and $J_0 \sim r^{-1/2} \cos(r - \frac{\pi}{4})$, for large $r$.

We will write the transformation as $u= \T_{\mathbb{R}^N} (v)$. Clearly, this is an isometry from the space \ $X=L^2(\mathbb{R}^N)$ \ into the space \ $\widetilde X=L^2(d\mu,\mathbb{R}^N))$, $d\mu=|x|^{2-N}\, J_0^2\, dx$. The great advantage of this formula is that it simplifies $I_{\mathbb{R}^N}(u)$, at least for smooth functions; For any $v \in C^{\infty}_0 (\mathbb{R}^N)$ holds that
\begin{eqnarray}\label{qf.vRN}
I_{\mathbb{R}^N} [u] &=& \int_{\mathbb{R}^N} |x|^{-(N-2)}\, J_0^2 (|x|)\, |\nabla v|^2\, dx \nonumber \\ &&+ \int_{\mathbb{R}^N} |x|^{-(N-2)}\, J_0^2 (|x|)\, v^2\, dx + \L(v),
\end{eqnarray}
where $\L$ is the Hardy singularity energy be defined in (\ref{Lfinite}). \vspace{0.3cm} \\
\emph{Proof of Lemma \ref{lemHPRN}}\ Let $v$ be a $C^{\infty}_0 (\mathbb{R}^N)$- function, and set
\begin{equation} \label{changevar}
u = |x|^{-\frac{N-2}{2}}\, J_0 (|x|)\, v,
\end{equation}
where $J_0$ is the Bessel function. If we do the change of variables (\ref{changevar}), we have that $J_0\, v \in C_0^{\infty}(\mathbb{R}^N)$, so from Subsection \ref{sub.definition}
\begin{equation}\label{J0v}
I_{\mathbb{R}^N} [u] = \int_{\mathbb{R}^N} |x|^{-(N-2)}\, \left ( \nabla \left (J_0 (|x|)\, v \right ) \right )^2\, dx + \L(v).
\end{equation}
After some calculations, we get that
\begin{eqnarray} \label{9.1}
I_{\mathbb{R}^N} [u] = \int_{\mathbb{R}^N} |x|^{-(N-2)}\, J_0^2 (|x|)\, |\nabla v|^2\, dx +  \int_{\mathbb{R}^N} |x|^{-(N-2)}\, (\nabla J_0 (|x|) )^2\, v^2\, dx \nonumber \\ + \int_{\mathbb{R}^N} |x|^{-(N-2)}\, J_0 (|x|)\, \nabla J_0 (|x|) \cdot \nabla v^2\, dx+ \L(v).
\end{eqnarray}
If we denote by $D_{\epsilon, R}$ to be the domain $\mathbb{R}^N \backslash \{B_{\epsilon} (0) \cup B_{R} (0)\}$, holds that
\begin{eqnarray*}
\int_{D_{\epsilon, R}} |x|^{-(N-2)}\, J_0 (|x|)\, \nabla J_0 (|x|) \cdot \nabla v^2\, dx &=& - \int_{D_{\epsilon, R}} \nabla \left ( |x|^{-(N-2)}\, J_0 (|x|)\, \nabla J_0 (|x|) \right )\, v^2\, dx \\ &&+ \int_{S_R} |x|^{-(N-2)}\, J_0 (|x|)\, \nabla J_0 (|x|)\, v^2\, \eta\, ds \\ &&- \int_{S_\epsilon} |x|^{-(N-2)}\, J_0 (|x|)\, \nabla J_0 (|x|)\, v^2\, \eta\,ds.
\end{eqnarray*}
Taking the limit as $R \to \infty$ we have that the second integral tends to zero, since $v \in C_0^{\infty}(\mathbb{R}^N)$. For the last integral, we have that $J_0' = -J_1$. Where $J_1$ satisfies $J_1 (0) =0$. Hence in the limit $\epsilon \downarrow 0$, the last integral also vanishes. Then, from (\ref{9.1}) we have that
\begin{eqnarray*}
I_{\mathbb{R}^N}[u] &=& \int_{\mathbb{R}^N} |x|^{-(N-2)}\, J_0^2 (|x|)\, |\nabla v|^2\, dx \nonumber \\ &&-  \int_{\mathbb{R}^N} |x|^{-(N-2)}\, J_0 (|x|)\, \left ( \Delta J_0 (|x|) - \frac{N-2}{|x|} J_0' (|x|) \right )\, v^2\, dx+ \L(v).
\end{eqnarray*}
However,
\[
\Delta J_0 (r) - \frac{N-2}{r} J_0' (r) = - J_0 (r),\;\;\; \mbox{for any}\;\; r.
\]
If we denote by $J_{\mathbb{R}^N}[v]$ the quantity
\[
J_{\mathbb{R}^N}[v] := \int_{\mathbb{R}^N} |x|^{-(N-2)}\, J_0^2 (|x|)\, |\nabla v|^2\, dx + \int_{\mathbb{R}^N} |x|^{-(N-2)}\, J_0^2 (|x|)\, v^2\, dx,
\]
we conclude that the Hardy functional is equal to
\begin{equation}\label{vhardyterm}
I_{\mathbb{R}^N} [u] = J_{\mathbb{R}^N}[v] + \L(v),
\end{equation}
thus
\begin{eqnarray*}
\int_{\mathbb{R}^N} |\nabla u|^2\, dx - \biggr ( \frac{N-2}{2} \biggr
)^2 \int_{\mathbb{R}^N} \frac{u^2}{|x|^2}\, dx &>& \int_{\mathbb{R}^N} |x|^{-(N-2)}\, J_0^2 (|x|)\, v^2\, dx \\
&=& \int_{\mathbb{R}^N} u^2\, dx.
\end{eqnarray*}
Finally, we have that
\begin{equation} \label{infRN}
\inf_{v \in C_0^{\infty}(\mathbb{R}^N)} \frac{\int_{\mathbb{R}^N} |x|^{-(N-2)}\, J_0^2 (|x|)\, |\nabla v|^2\, dx}{\int_{\mathbb{R}^N} |x|^{-(N-2)}\, J_0^2 (|x|)\, v^2\, dx} = 0.
\end{equation}
To see this we use the standard procedure; Let $v_n \in C_0^1 (\mathbb{R}^N)$ be such that $v_n (r) = C_1$ for $0<r<\frac{n\pi}{4}$, $|v_n'| < C_3$ for $\frac{n\pi}{4} <r<\frac{(n+1)\pi}{4}$ and $v_n(r)=0$ for $r>\frac{(n+1)\pi}{4}$. Then, using the fact that $J_0 (r) \sim r^{-1/2} \cos(r - \frac{\pi}{4})$, for $r>>1$, we have that
\[
\int_{\mathbb{R}^N} |x|^{-(N-2)}\, J_0^2 (|x|)\, |\nabla v_n|^2\, dx \leq c\, \int_{\frac{n\pi}{4}}^{\frac{(n+1)\pi}{4}} \cos^2 (r - \frac{\pi}{4})\, dr = C,
\]
for any $n$, large enough. However,
\[
\int_{\mathbb{R}^N} |x|^{-(N-2)}\, J_0^2 (|x|)\, v^2\, dx > c\, \int_{0}^{\frac{(n+1)\pi}{4}} \cos^2 (r - \frac{\pi}{4})\, dr \to \infty,
\]
as $n \to \infty$. Hence, $1$ is the best constant in (\ref{HardyPoincare}) and there exists no minimizer. $\blacksquare$

\begin{remark} \label{remR2}
The above arguments may be applied to the case of $D^{1,2}(\mathbb{R}^2)$: Let $v$ be any function in $C^{\infty}_0 (\mathbb{R}^2)$ and
\begin{equation} \label{transR2}
u = J_0 (|x|)\, v,
\end{equation}
Then,
\begin{equation} \label{PoincareR2}
\int_{\mathbb{R}^2} |\nabla u|^2\, dx > \int_{\mathbb{R}^2} u^2\, dx.
\end{equation}
The best constant in {\rm (\ref{PoincareR2})} is $1$ and there exists no minimizer.
\end{remark}
\begin{remark} \label{remnothing}
In order to clarify the transformation {\rm (\ref{transRN})} we argue as follows. Let $\Omega \subset \mathbb{R}^N$ be a bounded smooth domain. Assume that $(\lambda_1, u_1)$ be the principal eigenpair of the Laplace operator on $\Omega$. Then, for every $v \in C^{\infty}_0 (\Omega)$ and $u = u_1\, v$, holds that
\[
||u||^2_{H_0^1 (\Omega)} = \int_{\Omega} u_1^2\, |\nabla v|^2\, dx + \lambda_1\, \int_{\Omega} u_1^2\, |v|^2\, dx.
\]
In our case, the Bessel function $J_0$ solves the following eigenvalue problem
\[
-\Delta u = \lambda\, u,\;\;\; x \in \mathbb{R}^2,
\]
with $\lambda=1$.
\end{remark}

The following remark, makes clear that the transformation (\ref{transRN}) is not valid for every $u \in C_0^{\infty}(\mathbb{R}^N)$, so that $J[v]$ to be finite. Thus the Hardy-Poincar\'{e} inequality (\ref{HardyPoincare}), in Lemma \ref{lemHPRN}, makes sense for a certain class of functions. These functions, as it will be clear in the sequel, behave at the origin like the functions of $H$, at $\partial B_{z_m}(0)$, the zeroes of the Bessel function, must vanish like $||x| - z_m|^{a}$, $a \geq 1/2$, and at infinity must tend to zero.

\begin{remark} \label{remunotinJ}
The change of variables (\ref{changevar}) is not proper for every $u \in C_0^{\infty}(\mathbb{R}^N)$. For example, if $\overline{supp}\{u(r)\} \in (0,z_{0,1})$, then $v = |x|^{\frac{N-2}{2}}\, J_0^{-1} (|x|)\, u$ has $J[v] < \infty$. While, for $\overline{supp}\{u(r)\} = [\varepsilon,z_1 + \delta]$, $0<\delta < z_2 -z_1$, such that $u(z_1) \ne 0$, we have that $J[v] = \infty$.
\end{remark}

\medskip

\textbf{Functional spaces} Next, we discuss some properties of the corresponding functional space. We focus on the behavior at infinity. We introduce the space $\widetilde {\J} = W^{1,2} (d\mu, \mathbb{R}^N)$, $d\mu = |x|^{2-N}\, J^2_0 (|x|)\, dx$, to be defined as the closure of $C_0^{\infty}(\mathbb{R}^N)$ with respect to the norm
\[
||v||^2_{\widetilde {\J}} := J_{\mathbb{R}^N}[v] = \int_{\mathbb{R}^N} |x|^{-(N-2)}\, J_0^2 (|x|)\, |\nabla v|^2\, dx + \int_{\mathbb{R}^N} |x|^{-(N-2)}\, J_0^2 (|x|)\, v^2\, dx.
\]
We can easily verify that this weighted space $\widetilde {\J}$ is reasonably defined (see, for instance \cite{ok}), since the weight $|x|^{-(N-2)}\, J_0^2 (|x|)$ belongs to $L^1_{loc}$, with $|x|^{N-2}\, J_0^{-2} (|x|)$ belonging also to $L^1_{loc}$. The latter is true since the zeroes of the Bessel functions are simple. Moreover, as we will see below, every function belonging in $\widetilde {\J}$ tends to zero at infinity. Thus, $\widetilde {\J}$ is a well defined space.

We examine the behavior of the elements of $\widetilde {\J}$ at infinity. Note that close at the origin, $\widetilde {\J}$ is equivalent to $\widetilde {\cal H}$. \vspace{0.2cm} \\
\emph{Radial Case} \ \ Let $A = (R,+\infty) \backslash \bigcup (z_m-\varepsilon, z_m+\varepsilon) \subset \mathbb{R}$, $\varepsilon$ small enough. Assume that $v =v(r)$ is a radial function in $\J$ and $R$ is large enough. Then,
\[
||v||^2_{\widetilde {\J} (B_R^c)} \sim \int_{R}^{\infty} \cos^2 (r-\pi/4)\, \left ( (v')^2 + v^2 \right )\, dr.
\]
This means that in $B_R^c \backslash \bigcup \left ( B_{z_m +\varepsilon} (0) / B_{z_m -\varepsilon} (0) \right )$ the norm of $v$ is equivalent to $||v||_{H^1(A)}$. From this we conclude that
\[
\lim_{r \to \infty} v(r) = 0,
\]
up to the behavior at $r=z_m$. Moreover,
\[
v \in L^q \left ( B_{\varepsilon}^c \backslash \bigcup \left ( B_{z_m +\varepsilon} (0) / B_{z_m -\varepsilon} (0) \right ) \right ),
\]
for any $2 \leq q \leq \infty$. \vspace{0.1cm} \\
\emph{Nonradial Case}\ \ Let $v$ be a zero mean value function in $\widetilde {\J}$ and $u=\T_{\mathbb{R}^N} (v)$. Then, (\ref{vhardyterm}) holds with $\L(u)=0$. Since $u$ is nonradial, such that $I_{\mathbb{R}^N} [u] < \infty$, we may obtain that $u \in D^{1,2} (\mathbb{R}^N)$, thus $u \in L^{\frac{2N}{N-2}} (\mathbb{R}^N)$. Moreover, from Lemma \ref{lemHPRN} we have that $u$ belongs also to $L^{2} (\mathbb{R}^N)$. Finally, we get that $u = o (|x|^{-N/2})$ and so $v = o(|x|^{-1/2})$, as $|x| \to \infty$. \vspace{0.1cm} \\
\emph{General Case} \ \ The previous discussion was made in order to highlight some further properties of $\widetilde {\J}$. We can argue directly from the definition of the norm of $\widetilde {\J}$. Finally,
\begin{lemma}
Assume that $v \in \widetilde {\J}$. Then,
\[
\lim_{|x| \to \infty} v(x) = 0,
\]
up to the behavior at $\partial B_{z_m} (0)$.
\end{lemma}

We translate these results to the original framework. $\J$ is
defined as the isometric space of $\widetilde {\J}$ under the transformation $\T_{\mathbb{R}^N}$
given by (\ref{transRN}). The above definition is equivalent to say
that $\J$ is defined as the completion of the set
\[
\left \{ u = |x|^{-\frac{N-2}{2}}\, J_0 (|x|)\, v,\;\;\; v \in C_{0}^{\infty}
(\mathbb{R}^N) \right \}= \T_{\mathbb{R}^N} (C_{0}^{\infty}(\mathbb{R}^N)),
\]
under the norm $N_{\mathbb{R}^N} (u) = \|u\|_{\J}$ defined by
\begin{eqnarray} \label{newJRN}
||u||^2_{\J} &=:& \int_{\mathbb{R}^N} |x|^{-(N-2)}\,
J_0^2\, |\nabla (|x|^{\frac{N-2}{2}}\, J_0^{-1}\, u)|^2 \, dx \nonumber \\
&& + \int_{\mathbb{R}^N} |x|^{-(N-2)}\,
J_0^2\, ||x|^{\frac{N-2}{2}}\, J_0^{-1}\, u|^2 \, dx.
\end{eqnarray}
Note that the last integral is nothing else than the $L^2 (\mathbb{R}^N)$-norm. Thus, for any $u \in \J$, we have that $u$ decays to zero at infinity. \vspace{0.3cm}

\textbf{The case of $\J$} We discuss the connection of $\J$ with the Hardy functional. We also give the properties of the functions belonging to $\J$.
\begin{lemma}
Let $B_{\varepsilon,m} = B_{\varepsilon} (0) \bigcup_{m=1} B_{z_m +\varepsilon} (0) / B_{z_m -\varepsilon} (0)$, $z_m$ are the zeroes of the Bessel function. For any $u \in \J$, holds that
\begin{equation} \label{normJHR}
||u||_{\J}^2 = \lim_{\varepsilon \to 0} I_{B^c_{\varepsilon,m}} [u] - \L_{\varepsilon}(u) + \sum_{m=1}^{\infty} \L_{\varepsilon^+,m}(u) - \L_{\varepsilon^-,m}(u),
\end{equation}
where, $B^c_{\varepsilon,m}$ denotes the complement of $B_{\varepsilon,m}$ in $\mathbb{R}^N$, $\L_{\varepsilon}$ is the Hardy singularity energy at the origin, defined in {\rm (\ref{Lepsilon})}, and
\[
\L_{\varepsilon^+,m}(u) = \int_{\partial B_{z_m + \varepsilon}(0)} J_0^{-1}\, J'_0\, u^2\, ds,\;\;\;\;\;\; \L_{\varepsilon^-,m}(u) = \int_{\partial B_{z_m - \varepsilon}(0)} J_0^{-1}\, J'_0\, u^2\, ds.
\]
The quantities $\L_{\varepsilon^+,m}(u)$ and $-\L_{\varepsilon^-,m}(u)$ are always positive.
\end{lemma}
\emph{Proof} The proof relies on the careful examination of the surface integrals arising in the calculations. Let $\varepsilon$ and $R$ be a small and a large, respectively, enough positive numbers. If $B_{\varepsilon, R} = B_R \backslash B_{\varepsilon,m} $. Calculating the norm of $u$ in $\J$ on $B_{\varepsilon, R}$ we face the existence of the following surface integrals:
\[
\int_{S_{\varepsilon}} |x|^{-1}\, u^2\, ds,\;\;\;\;\;\; \int_{S_{R}} |x|^{-1}\, u^2\, ds,
\]
and
\[
\int_{\partial B_{\varepsilon,m}} J_0^{-1}\, J'_0\, u^2\, ds,\;\;\;\;\;\; \int_{S_{R}} J_0^{-1}\, J'_0\, u^2\, ds.
\]
For $\varepsilon$ small enough, the first surface integral corresponds to $\L_{\varepsilon}$. The third corresponds to $\L_{\varepsilon^+,m}(u)$ and $\L_{\varepsilon^-,m}(u)$. Observe that $J'_0 = -J_1$ and that $J_0$ and $J_1$ have the opposite sign at $z_m +\varepsilon$, while they have the same sign at $z_m - \varepsilon$; thus, the quantities $\L_{\varepsilon^+,m}(u)$ and $-\L_{\varepsilon^-,m}(u)$ are always positive.

For $R$ large enough, the key observation is that the definition of the space $\J$ implies that $u \in L^2 (\mathbb{R}^N)$; the second and the fourth surface integrals tend to zero. Finally, we conclude that (\ref{normJHR}) is true. $\blacksquare$ \vspace{0.2cm}

Next, we investigate some further properties of the functions that belong in $\widetilde {\J}$ and $\J$. Observe that if $v \in \widetilde {\J}$ and $u = \T_{\mathbb{R}^N} (v)$, from (\ref{newJRN}) we have that $||v||_{\widetilde {\J}} = ||u||_{\J} \geq ||u||_{L^2 (\mathbb{R}^N)}$. Let $D_{\varepsilon,m} = B_{z_{m+1} -\varepsilon} (0) / B_{z_m +\varepsilon} (0)$, for some fixed $m$. Then,
\begin{equation} \label{Iem}
I_{D_{\varepsilon,m}}[u] > - c_* \int_{D_{\varepsilon,m}} \frac{u^2}{|x|^2}\, dx \geq - c\, ||u||^2_{L^2 (D_{\varepsilon,m})},
\end{equation}
i.e. is bounded from below, for every $m$. Assume, for the moment, that $u \equiv 0$, in $B_{\varepsilon} (0)$. From (\ref{Iem}) and (\ref{normJHR}) we conclude that $u \in \J$, or equivalently $v \in \widetilde {\J}$, if and only if
\begin{equation} \label{Lem}
\sum_{m=1}^{\infty} \L_{\varepsilon^+,m}(u) - \L_{\varepsilon^-,m}(u) < \infty.
\end{equation}
This actually means that
\begin{equation} \label{uzm0}
\lim_{\varepsilon \to 0} u(z_m + \varepsilon) = \lim_{\varepsilon \to 0} u(z_{m+1} - \varepsilon) = 0,
\end{equation}
for any $m$, and particular
\begin{equation} \label{uz1m0}
\lim_{\varepsilon \to 0} u(z_1 - \varepsilon) =0.
\end{equation}
Then $u \in \H_0^1 (D_{\varepsilon,m})$, for any $m$, and $I_{B^c_{\varepsilon} (0)} [u]$ is a well defined positive number.

The behavior at the origin, as it is expected, is a local effect and is described in Section 2. To see this we may use the standard argument using cutoff functions; $I_{B_{\varepsilon} (0)} - \L_{\varepsilon}$ is at least bounded from below. We finally have
\begin{lemma} \label{lemmabehaviorHR}
For every $u \in \J$ and every $\varepsilon,\, m$, the quantity $I_{B^c_{\varepsilon,m}} [u]$ is well defined and positive. Then, $u \in \J$ if and only if (\ref{Lem}) holds. Moreover, $u$ vanishes at $z_m$ with rate $O(||x|-z_m|^a)$, $a \geq 1/2$. In particular, the restriction of $u$ on $B_{z_1}(0)$ is a function that belongs to $H$, i.e.
\[
\J (B_{z_1} (0))\;\;\; \mbox{is equivalent to}\;\;\; H (B_{z_1} (0)).
\]
Its behavior at the origin was studied in Section 2. Finally, $u \in L^2 (\mathbb{R}^N)$, so it tends to zero at infinity. As it follows,
\begin{equation} \label{uH1R}
\mbox{if}\;\;\; u \in \J\;\;\; \mbox{then}\;\;\; u \in H^1 (B^c_{z_1}),
\end{equation}
i.e. $\J(B^c_{z_1}) \subset H^1 (B^c_{z_1})$.
\end{lemma}

\medskip

\textbf{Compact imbeddings}\ \ From the above discussion and the compact results of (\ref{compLp}), we have the following compact imbeddings
\begin{itemize}
  \item $\J_r(B^c_{z_1}) \hookrightarrow L^q_r (B^c_{z_1})$, $2<q<\frac{2N}{N-2}$, where $J_r$ and $L^q_r$ are the radial subspaces of $J_r$ and $L^q_r$, respectively.
  \item $\J (B_{z_1} (0)) \hookrightarrow L^p (\Omega)$, $1\leq p < \frac{2N}{N-2}$,
  \item $\J \hookrightarrow L^q_g (\mathbb{R}^N)$, $1 \leq q < \frac{2N}{N-2}$, where $g$ is a weight function belonging in $L^{N/2} (\mathbb{R}^N)$.
\end{itemize}

\medskip

\textbf{Conclusion}\ \ The space $\J$ is a well defined space, through $\J$, as the isometric space of $\widetilde {\J}$ under the transformation $\T_{\mathbb{R}^N}$. In this space there exist ``bad functions'', concerning the behavior at the origin, as was discussed in Section \ref{sec.pfs}. Concerning the behavior at infinity, each function belongs in $L^2(\mathbb{R}^N)$ and in $H^1 (B^c_{z_1})$. In addition, if $u \in \J$, must vanish at the zeroes of the Bessel function, as Lemma \ref{lemmabehaviorHR} states.

As we can see, the transformation (\ref{transRN}) removes the singularity at infinity. However, the cost we have to pay, is the existence of an infinite series of singularities at $z_m$. Contrary to the case of the origin, there exists a singularity energy at these points, which is possible to increase the $I[u]$-"energy".

The study of the Cauchy problem (\ref{eq1}) on $\mathbb{R}^N$ is not in the purpose of this work; actually this was done in \cite{vz00}. However, we make some comments. This problem is well defined on $\J$. An important observation is that when $u(x,t)$ is a solution of equation (\ref{eq1}), then $v$ satisfies the following associated equation
\begin{equation} \label{Cauchy_v}
v_t = |x|^{N-2}\; J_0^{-2}\; \nabla \cdot \left( |x|^{-(N-2)}\, J_0^2\, \nabla v \right) -v,
\end{equation}
with clear equivalence for $x \ne 0, z_m$. The associated space is $\widetilde {\J}$, with base space
$\widetilde X_{\mathbb{R}^N}=L^2(d\mu,\mathbb{R}^N)$, $d\mu = |x|^{-(N-2)}\, J_0^2$, with
\[
\widetilde {\J} \subset \widetilde X_{\mathbb{R}^N} \subset \widetilde {\J}'.
\]
However, as it seems, the above problem may be not appropriate in order to study (\ref{eq1}). First, the positive solutions of (\ref{Cauchy_v}) correspond, through (\ref{transRN}), to solutions of (\ref{eq1}) which change sign. In addition, these solutions must vanish at $z_m$, thus, the dynamic of (\ref{eq1}) described by (\ref{Cauchy_v}) may be nothing else but the dynamics of (\ref{eq1}) defined on each bounded domain $B_{z_m+1} (0) \backslash B_{z_m} (0)$, up to a scaling.

\section{Appendix} \label{sec.app}
\setcounter{equation}{0}
%
%
\subsection {Alternative proof of Proposition \ref{propH=H}}

We will give the proof of Proposition \ref{propH=H} using a different approach from this of Subsection \ref{sub.H=H}. \vspace{0.2cm} \\
\emph{Alternative Proof of Lemma \ref{lemmacalHcontains}} Let $v$ be a function such that $v|_{\partial \Omega} =0$ and $||v||_{\widetilde {\cal H}}$ holds.We will find a $C_0^{\infty} (\Omega)$-function such that
\begin{equation} \label{H1.2}
\int_{\Omega} |x|^{-(N-2)}\, |\nabla \phi - \nabla v|^2\, dx < \varepsilon,
\end{equation}
for any $\varepsilon>0$. Note that since $||v||_{\widetilde {\cal H}}$ holds, we get that $v \in H_0^1 (\Omega)$, so $v \in L^2 (\Omega)$. Thus $v$ may be written as the sum of its radial part $v_r$ and the nonradial part $v_{nr}$. If we denote the same for $\phi$ we have that
\begin{eqnarray} \label{H1.3}
\left ( \int_{\Omega} |x|^{-(N-2)}\, |\nabla \phi - \nabla v|^2\, dx \right )^{1/2} &\leq& \left ( \int_{\Omega} |x|^{-(N-2)}\, |\nabla \phi_r - \nabla v_r|^2\, dx \right )^{1/2} \nonumber \\ && + \left ( \int_{\Omega} |x|^{-(N-2)}\, |\nabla \phi_{nr} - \nabla v_{nr}|^2\, dx \right )^{1/2}.
\end{eqnarray}
The first integral in (\ref{H1.3}) is equal (up to some constant) to
\[
\int_{B_R \subset \mathbb{R}^2} |\nabla \phi - \nabla v|^2\, dx,
\]
for some $R>0$, large enough. Since the $C_0^{\infty} (B_R)$ functions are dense in $H_0^1 (B_R \subset \mathbb{R}^2)$ we have the existence of a $C_0^{\infty} (B_R)$ function such that
\[
\left ( \int_{B_r} |x|^{-(N-2)}\, |\nabla \phi_r - \nabla v_r|^2\, dx \right )^{1/2} \leq \frac{1}{2}\, \varepsilon^{1/2},
\]
for any $\varepsilon >0$. It is clear that this function is a mollifier of $v_r$ in $B_R \subset \mathbb{R}^2$ and may be chosen to be radial belonging in $C_0^{\infty} (\Omega)$. For the second integral in (\ref{H1.3}) we have that $\phi_{nr} - v_{nr}$ has zero radial part. Then, $|x|^{-(N-2)/2} (\phi_{nr} - \nabla v_{nr})$ belongs in $H^1 (\Omega)$. More precisely,
\begin{eqnarray*}
\int_{\Omega} |x|^{-(N-2)}\, |\nabla \phi_{nr} - \nabla v_{nr}|^2\, dx &=& I[|x|^{-(N-2)/2} (\phi_{nr} - v_{nr})] \\ &\leq& c\, |||x|^{-(N-2)/2} (\phi_{nr} - v_{nr})||^2_{H^1_0 (\Omega)},
\end{eqnarray*}
where $I$ is the Hardy functional. However, the $C_0^{\infty} (\Omega \backslash \{0\})$- functions are dense in $H_0^1 (\Omega)$, for instance see (\ref{H1}). Thus there exists a function $\phi_{nr}$, belonging in $C_0^{\infty} (\Omega \backslash \{0\})$, such that
\[
\left ( \int_{\Omega} |x|^{-(N-2)}\, |\nabla \phi_{nr} - \nabla v_{nr}|^2\, dx \right )^{1/2} \leq \frac{1}{2}\, \varepsilon^{1/2},
\]
for any $\varepsilon >0$. Is clear from the construction of Lemma \ref{denseH} that $\phi_{nr}$ has zero mean value. Finally, setting $\phi = \phi_r + \phi_{nr}$, we conclude that (\ref{H1.2}) holds for any $\varepsilon >0$ and the proof is completed. $\blacksquare$ \vspace{0.2cm} \\
We define the space $\hat{\H}$ as the completion of the set
\[
B:= \left \{ w ;\;\;\; w = |x|^{\frac{N-2}{2}}\, \phi, \;\;\; \phi \in C_0^{\infty}(\Omega)  \right \}
\]
under the norm
\[
||w||^2_{\hat{\H}} = \int_{\Omega} |x|^{-(N-2)}\, \left | \nabla \left ( |x|^{\frac{N-2}{2}}\, \phi \right ) \right |^2\, dx\,.
\]
From the fact that
\[
||w||^2_{\hat{\H}} = I_{\Omega}(\phi),
\]
for any $\phi \in C_0^{\infty}(\Omega)$, we obtain that this space is isometric with $H$.

Next we give the relation between these two spaces: $\widetilde {\cal H}$ and $\hat{\H}$. Note that for every function in the set $B$, we have that $||w||_{\hat{\H}} = ||w||_{\widetilde {\cal H}}$. Thus, from
Lemma \ref{lemmacalHcontains} $w$ belongs in $\widetilde {\cal H}$ and since the two norms coincide we have that
\[
\hat{\H} \subseteq \widetilde {\cal H}
\]
Moreover, we have the following
\begin{lemma} \label{hatH=H1}
The space $\hat{\H}$ is not a strict subspace of $\widetilde {\cal H}$.
\end{lemma}
\emph{Proof} Assume the opposite. Then there exists function $v
\in \widetilde {\cal H}$, $v \not\equiv 0$, such that
\[
v \perp \hat{\H},
\]
which means that
\[
<v,w>_{\widetilde {\cal H}} = 0,\;\;\; \mbox{for any}\;\;\; w \in
\hat{\H}.
\]
The latter is equivalent to
\[
\int_{\Omega} |x|^{-(N-2)}\, \nabla v \cdot \nabla \left ( |x|^{\frac{N-2}{2}}\, \phi \right )\, dx = 0,\;\;\; \mbox{for any}\;\;\; \phi \in C_0^{\infty}(\Omega).
\]
Integrating by parts and noting that the boundary terms are zero, we have that
\[
\int_{\Omega} div \left (|x|^{-(N-2)}\, \nabla v \right ) |x|^{\frac{N-2}{2}}\, \phi \, dx = 0,\;\;\; \mbox{for any}\;\;\; \phi \in C_0^{\infty}(\Omega),
\]
or
\begin{equation}\label{eq2}
\int_{\Omega} |x|^{-\frac{N-2}{2}}\, \left ( |x|^{(N-2)}\, div \left (|x|^{-(N-2)}\, \nabla v \right ) \right )\, \phi \, dx = 0,\;\;\; \mbox{for any}\;\;\; \phi \in C_0^{\infty}(\Omega).
\end{equation}
However, we know that the operator $L$ with domain
\[
D_* = \left \{ f \in \widetilde {\cal H}: |x|^{N-2}\; \nabla \cdot
\left( |x|^{-(N-2)}\, \nabla f \right) \in L^2(|x|^{2-N}dx,\Omega)
\right \}.
\]
is self-adjoint and bijective with compact inverse. Thus
\[
|x|^{-\frac{N-2}{2}}\, \left ( |x|^{(N-2)}\, div \left (|x|^{-(N-2)}\, \nabla v \right ) \right ) \in L^2 (\Omega)
\]
and from (\ref{eq2}) we conclude that
\[
|x|^{-\frac{N-2}{2}}\, \left ( |x|^{(N-2)}\, div \left (|x|^{-(N-2)}\, \nabla v \right ) \right ) = 0,
\]
almost everywhere in $\Omega$. The fact that $L$ is bijective implies that $v \equiv 0$ which is a contradiction and the proof is completed. $\blacksquare$ \vspace{0.2cm}

As a conclusion we have that $\hat{\H}$ coincides with $\widetilde {\cal H}$, i.e., $\hat{\H} = \widetilde {\cal H}$. This in terms of $u$, i.e. in terms of the spaces $\H$ and $H$, means that these two spaces coincide.

\vspace{1cm}

\noindent \textsc{Acknowledgment.}  The first author supported by Spanish
Project MTM2008-06326-C02-01. Work partly done during a visit of the second author to Univ. Au´\'onoma de Madrid,
supported by the same project.

\

%

{\small
\bibliographystyle{amsplain}

\begin{thebibliography}{10}
%
%
\bibitem{BG} P. Baras, J. A. Goldstein, \emph{The heat equation with a singular potential}.
Trans. Amer. Math. Soc.  {\bf 284}  (1984),  no. 1, 121--139.
%
\bibitem{bft} G. Barbatis, S. Filippas and A. Tertikas \emph{Critical heat kernel estimates for Schr\"{o}}dinger operators via Hardy-Sobolev inequalities, J. Funct. Analysis  \textbf{208} (2004), 1-30.
%
\bibitem{bm} H. Brezis and M. Marcus, Hardy's inequality revisited,
{\em Ann. Sc. Norm. Pisa} {\bf 25} (1997), 217-237.
%
\bibitem{bv97} H. Brezis and J. L. V\'{a}zquez,
Blowup solutions of some nonlinear elliptic problems, {\em Revista
Mat. Univ. Complutense Madrid} \textbf{10} (1997), 443-469.
%
\bibitem{cm99} X. Cabr\'{e}, Y. Martel \emph{Existence versus explosion instantan\'{e} pour des equations de lachaleur lin\'{e}aires avec potentiel singulier}, C.R. Acad. Sci. Paris \textbf{329} (1999), 973-978.
%
\bibitem{ckn} L. Caffarelli, R. Kohn and L. Nirenberg, \emph{First
Order Interpolation Inequalities with Weights}, Compositio Math.
\textbf{53} (1984), no 3., 259-275.
%
\bibitem{cw01} F. Catrina and Z.-Q. Wang, \emph{On the
Caffarelli-Kohn-Nirenberg Inequalities: Sharp Constants, Existence
(and Nonexistence) and Symmetry of Extremal Functions}, Comm. Pure
Appl. Math. \textbf{LIV} (2001), 229-258.
%
\bibitem{cf08} A. Cianchi and A. Ferone, \emph{Hardy inequalities with non-standard remainder terms}, Ann. I. H. Poincar\'{e} – AN  \textbf{25} (2008), 889-906.
%
\bibitem{d99} E. B. Davies, {\em A review of Hardy inequalities,}  Oper. Theory Adv.
Appl. {\bf 110} (1999), 55-67.
%
\bibitem{kmp} A. Kufner, L. Maligranda and L.-E. Persson,
{\em The Hardy inequality. About its history and some
related results}, Vydavatelsky' Servis, Plzen, 2007.
%
%
\bibitem{maz85} V. G. Maz'ja, \emph{``Sobolev spaces''}, Springer-Verlag, 1985.
%
\bibitem{ok} B. Opic, and A. Kufner, Hardy type inequalities,
Pitman Rechearch Notes in Math., {\bf 219} Longman 1990.
%
\bibitem{vzog} J. L. V\'azquez and N. B. Zographopoulos, \emph{Functional
aspects of Hardy type inequalities}, in preparation.
%
\bibitem{vz00} J. L. V\'azquez and E. Zuazua, \emph{The Hardy
inequality and the asymptotic behaviour of the heat equation with
an inverse-Square Potential},  J. Functional Analysis \textbf{173}
(2000), 103-153.
%
\bibitem{z10} N. B. Zographopoulos, \emph{Existence of Extremal Functions for a Hardy-Sobolev
Inequality}, J. Funct. Analysis \textbf{259} (2010), 308-314.
%
%
\end{thebibliography}

}

\vskip 0.5cm

{\sc Addresses:}

{\sc Juan Luis V{\'a}zquez}\newline
Departamento de Matem\'{a}ticas, Universidad Aut\'{o}noma de Madrid, 28049
Madrid, Spain.  \newline Second affiliation: Institute ICMAT. \newline
e-mail: juanluis.vazquez@uam.es

\medskip

{\sc N. B. Zographopoulos}\newline
University of Military Education, Hellenic Army Academy, Department of Mathematics \& Engineering Sciences, Vari - 16673, Athens Greece,
\newline Second affiliation: Division of Mathematics, Technical University of Crete, 73100, Chania Greece
\newline e-mail: nzograp@gmail.com, zographopoulosn@sse.gr

\

{\bf Keywords:} \ {Hardy inequality, Inverse-Square Potential, Energy at singularity, Hidden or Dark energy.} \vspace{0.1cm}

\medskip

{\bf AMS Subject Classification (2010):} \  {???.}
\end{document}